\documentclass[10pt]{article}
\usepackage[utf8]{inputenc}
\usepackage[top=1in, bottom=1in, left=1in, right=1in]{geometry}
\usepackage[title,titletoc,toc]{appendix}
\usepackage{amssymb}
\usepackage{amsmath}
\usepackage{color}
\usepackage{tikz}
\usepackage{amsthm}
\newtheorem*{rk}{Remark}
\newtheorem{thm}{Theorem}

\newtheorem{prop}{Proposition}

\newtheorem{lemme}{Lemma}

\title{A finite speed of propagation approximation for the incompressible Navier-Stokes equations}
\author{Im\`ene HACHICHA \footnote{Laboratoire Analyse et Probabilités, Université d'Evry-Val d'Essonne 23 bd de France, 91037 Evry Cedex. E-mail : imene.hachicha@univ-evry.fr}}
\date{}

\newcommand{\aq}{A_{\mathbb{Q}}}
\newcommand{\bq}{B_{\mathbb{Q}}}
\newcommand{\ap}{A_{\mathbb{P}}}
\newcommand{\bp}{B_{\mathbb{P}}}
\newcommand{\q}{\mathbb{Q}}
\newcommand{\p}{\mathbb{P}}
\renewcommand{\u}{u^{\varepsilon , \alpha}}
\newcommand{\w}{w^{\varepsilon , \alpha}}
\newcommand{\z}{z^{\varepsilon , \alpha}}
\newcommand{\eps}{\varepsilon}
\renewcommand{\div}{\textrm{div}\, }
\newcommand{\R}{\mathbb{R}}
\newcommand{\dt}{\partial_{t}}
\newcommand{\dtt}{\partial_{tt}}
\renewcommand{\d}{\textrm{d}}

\newcommand{\epsalpha}{{\eps , \alpha}}

\begin{document}

\maketitle
\begin{abstract}
In this paper, we introduce a finite propagation speed perturbation of the incompressible  Navier-Stokes equations $(NS)$. The model we consider is inspired by a hyperbolic perturbation of the heat equation due to Cattaneo in \cite{cattaneo} and by an equation that Vi\v{s}ik and Fursikov investigated in \cite{vishik} in order to find statistical solutions to $(NS)$. We prove that the solutions to the perturbed Navier-Stokes equation approximate those to $(NS)$. 

We use refined energy methods involving fractional Sobolev spaces and precise estimates on the nonlinear term due to the dyadic Littlewood-Paley decomposition. 
\end{abstract}
\tableofcontents
\section{Introduction}
The purpose of this paper is to approximate the solutions to the incompressible Navier-Stokes equations with quasi-critical regularity initial datum by solutions to a nonlinear wave equation with a finite speed of propagation which is obtained by penalizing the incompressibility constraint.
First, let us recall the Navier-Stokes equations which govern the motion of an incompressible, viscous and homogeneous newtonian fluid whose velocity and pressure are denoted by $v$ and $p$ respectively. 
$$ (NS) ~~~ \partial_t v(t,x) - \nu \Delta v(t,x) + (v.\nabla)v(t,x) = -\nabla p(t,x) ~, ~~ \div v(t,x) =0 ,
$$
where $t>0$ and $x \in \R^n$ for $n=2,3$. The velocity is a $\R^n$ valued vector field and the pressure is scalar. The coefficient of the Laplacian is the viscosity and, without loss of generality, is assumed to be $1$ in the following.\\%by a scale change in the time variable, is assumed to be $1$.\\
%\red{(Under standard assumptions on $p$)} 
Applying the Leray projector $\p$ which maps $L^2$  into $L^2_\sigma:= \left\lbrace u \in L^2 : \div u =0 \right\rbrace$ to $(NS)$, we obtain the equations
$$ (NS) ~~~~~ \partial_t v - \Delta v + \p (v.\nabla)v = 0 ~, ~~ \div v =0
$$
from which we can recover the pressure $p$.\\
A first hyperbolic perturbation of $(NS)$ has been obtained after relaxation of the Euler equations and rescaling variables (see \cite{brenier} and references therein):
$$ (HNS^\eps) ~~~ \eps \partial_{tt}u^\eps + \partial_t u^\eps - \Delta u^\eps + \p (u^\eps.\nabla)u^\eps = 0 ~, ~~ \div u^\eps =0 .
$$
In \cite{cattaneo}, Cattaneo introduced this equation (without the nonlinear term) as a perturbation of the heat equation. In \cite{brenier} and \cite{paicu}, the authors approximate the solutions to $(NS)$ by solutions to $ (HNS^\eps)$ under some assumptions on the size and the regularity of the initial data. In \cite{article}, we improve the results of \cite{brenier} and \cite{paicu}. Under weaker assumptions on the initial data size, we prove the convergence of solutions to $ (HNS^\eps)$ towards solutions to $(NS)$ in the critical Sobolev space norm, that is $\dot H^{\frac{n}{2}-1}(\R^n)$, as $\eps$ goes to $0$. More precisely, the theorem is:
\begin{thm} \label{th}
 Let $n=2$ or $3$ and $0< s,\delta <1$. Let $v_0 \in H^{\frac{n}{2}-1+s} (\R^n)^n$ be a divergence-free vector field and $(u_0^\eps , u_1^\eps) \in H^{\frac{n}{2}+\delta}(\R^n)^n \times H^{\frac{n}{2}-1+\delta}(\R^n)^n$ be a sequence of initial data for problem $(HNS^\eps)$. Assume
$$
\left\{\begin{array}{ccc} 
\!\!\! \|u_0^\eps -v_0\|_{\dot H^{\frac{n}{2}-1}} + \eps \|u_1^\eps\|_{\dot H^{\frac{n}{2}-1}} + \eps^{\frac{1}{2}} \|u_0^\eps\|_{\dot H^{\frac{n}{2}}} + \eps^{\frac{1+\delta}{2}} \|u_0^\eps \|_{\dot H^{\frac{n}{2}+\delta}} + \eps^{\frac{\delta}{2}} \|u_0^\eps\|_{\dot H^{\frac{n}{2}-1+\delta}}= \mathcal{O} \left(\eps^{\frac{s}{2}}\right) \\
 \eps^{1+\frac{\delta}{2}} \|u_1^\eps \|_{\dot H^{\frac{n}{2}-1+\delta}} = o \left(1\right).
\end{array}\right.
$$
Moreover, if $n=3$, assume that $\|u_0^\eps \|_{\dot H^{\frac{1}{2}}(\R^3)^3} < \frac{1}{16}$.\\
Then, for $\eps$ small enough, there exists a global solution $u^\eps$ to system $(HNS^\eps)$ that converges, when $\eps$ goes to $0$, in the  $L_{\textrm{loc}}^\infty(\R^+; \dot H^{\frac{n}{2}-1}(\R^n)^n) $ norm, towards the unique solution $v$ to $(NS)$, with $v_0$ as initial datum. Moreover, for all positive $T$, there exists a constant $C_T$, depending only on $T$ and $v$,  such that
$$
\sup_{t\in [0,T]}  \|u^\eps -v\|_{\dot H^{\frac{n}{2}-1}(\R^n)^n}^2 \, dx \leq C_T \eps^{\left( \frac{s}{2} \right)^-}.
$$
\end{thm}
\begin{rk}
As a consequence of the assumptions $\|u_0^\eps \|_{\dot H^{\frac{1}{2}}(\R^3)} < \frac{1}{16}$ and $\|u_0^\eps -v_0\|_{ \dot H^{\frac{1}{2}}(\R^3) }  = \mathcal{O} \left(\eps^{\frac{s}{2}}\right) $, we obtain the smallness of $\|v_0 \|_{\dot H^{\frac{1}{2}}}$, which is a necessary condition to the existence of a global strong solution to the Navier-Stokes equations in $\R^3$.
\end{rk}
In this paper, we shall consider another hyperbolic perturbation of $(NS)$ which is inspired by Cattaneo's perturbation $(HNS^\eps)$ and by Vi\v{s}ik and Fursikov's weakly compressible equations in \cite{vishik}. The point is that we try to get around the difficulties which come from the Leray projector $\p$ (or, equivalently, from the pressure) on the one hand and, on the other hand, from the infinite propagation speed due to the heat kernel. The equation we introduce in this paper is:
$$
(HNS^\epsalpha) ~~~ \eps \partial_{tt} u^{\eps , \alpha} + \partial_t u^{\eps , \alpha} - \Delta u^{\eps , \alpha} = -(\u . \nabla) \u + \frac{1}{\alpha} \nabla (\div u^{\eps , \alpha}).
$$
Under the same assumptions as in Theorem \ref{th}, we prove that the solutions to $(HNS^\epsalpha)$ approximate those to $(NS)$. Notice that we do not need any more restrictions involving $\alpha$ on the initial data.
\begin{thm} \label{th_alpha}
Let $n=2$ or $3$ and $0< s,\delta <1$. Let $v_0 \in H^{\frac{n}{2}-1+s} (\R^n)^n$ be a divergence-free vector field and $(\u_0 , \u_1)=(u_0^\eps ,u_1^\eps) \in H^{\frac{n}{2}+\delta}(\R^n)^n \times H^{\frac{n}{2}-1+\delta}(\R^n)^n$ be a sequence of divergence-free initial data for problem $(HNS^\epsalpha)$, independent of $\alpha$. Assume
\begin{equation} \label{size}
\left\{\begin{array}{rcc} 
\|\u_0  -v_0\|_{\dot H^{\frac{n}{2}-1}} + \eps \|\u_1\|_{\dot H^{\frac{n}{2}-1}} + \eps^{\frac{1}{2}} \|\u_0\|_{\dot H^{\frac{n}{2}}} & = & \mathcal{O} \left(\eps^{\frac{s}{2}}\right) \\
\eps^{\frac{1+\delta}{2}} \|\u_0 \|_{\dot H^{\frac{n}{2}+\delta}} + \eps^{\frac{\delta}{2}} \|\u_0\|_{\dot H^{\frac{n}{2}-1+\delta}} & = & \mathcal{O} \left(\eps^{\frac{s}{2}}\right) \\
 \eps^{1+\frac{\delta}{2}} \|\u_1 \|_{\dot H^{\frac{n}{2}-1+\delta}} &=& o \left(1\right).
\end{array}\right.
\end{equation}
Moreover, if $n=3$, assume that $\|\u_0 \|_{\dot H^{\frac{1}{2}}(\R^3)^3} < \frac{1}{36 K_2^3}$, where $K_2$ is the constant such that 
\[
 \|f\|_{L^3(\R^3)} \leq K_2 \|f\|_{\dot H^{\frac{1}{2}}(\R^3)}.
\]
\\
Then, for $\eps , \alpha$ small enough and for all positive $T$, the global solutions to $(HNS^\epsalpha)$ approximate those to $(NS)$ in the  $L^\infty_T  \dot H^{\frac{n}{2}-1}(\R^n)^n $ norm.
%\red{ \begin{equation}
%  \label{size}
%  \sqrt{\alpha} \|\u_0\|_{L^2(\R^2)^2} \leq \frac{1}{2C},
% \end{equation}
%2C is the constant of Gagliardo-Nirenberg.}
\end{thm}
\begin{rk}
 A possible choice of initial data is $\widehat{u^{\epsalpha}_0}(\xi) = \widehat{ v_0} (\xi) \mathbf{1}_{\{\sqrt{\eps} |\xi|<1\}}$ and $\u_1 =0$.
\end{rk}

In this paper, we shall prove that the solutions to $(HNS^\epsalpha)$ with initial data $(u^\eps_0 , u^\eps_1)$ converge, as $\alpha$ goes to $0$, towards the solutions to $(HNS^\eps)$ with the same initial data. Then we conclude according to Theorem \ref{th}.\\

The next section is dedicated to the introduction of the model. In section \ref{vitessefinie}, we shall prove that $(HNS^\epsalpha)$ has a finite speed of propagation using the results of section \ref{existencelocale} where we prove local existence. Then,we focus on the 2D case in section \ref{2d}: first, we recall in subsection \ref{rappels2d} some important estimates on the solutions $u^\eps$ to $(HNS^\eps)$ then we prove global existence for $(HNS^\epsalpha)$ in part \ref{glob2d} and we show that its solutions approximate those to $(HNS^\eps)$ in subsection \ref{conv2d}. Finally, section \ref{3d} is devoted to the 3D case and follows the same plan as section \ref{2d}: in subsection \ref{rappels3d} we recall the useful regularity results on $u^\eps$ then we prove that the local solutions to $(HNS^\epsalpha)$ are global in subsection \ref{glob3d} and that the global solutions approximate those to $(HNS^\eps)$ in subsection \ref{conv3d} . Finally, Theorem \ref{th} allows to conclude the proof of Theorem 
\ref{th_alpha}. Some important estimates coming from the framework of Littlewood-Paley theory are recalled in appendix \ref{lptheory}.

\section{Introducing the model } \label{modele}
In this section we shall introduce the finite speed of propagation equation that we will study in the next sections. We will work in the setting of $\R^2$. The 3D case is similar.\\
First, let us perturb the Navier-Stokes system $(NS)$ into the damped nonlinear wave equation $ (HNS^\eps)$ which we recall:
$$ (HNS^\eps) ~~~
\eps \partial_{tt}u^\eps + \partial _t u^\eps - \Delta u^\eps = - (u^\eps.\nabla)u^\eps - \nabla p^\eps ~~,~~~ \div u^\eps = 0. 
$$
Then, applying the $\div =\sum_{i=1}^n \partial_i (.)_i$ operator to this equation , we obtain 
\begin{equation} \label{pf} 0 = \div f(u^\eps) - \Delta p^\eps , \end{equation} where $f(u)=-(u.\nabla)u$. 
Let us now consider the stationary problem
\begin{equation} \label{stat}
- \Delta w^\eps = f - \nabla p^\eps ~~,~~~ \div w^\eps = 0
\end{equation}
which is related to \eqref{pf}. Indeed, applying $\div$ to \eqref{stat}, we obtain \eqref{pf}. Given $f\in \dot H^{-1}$, we look for a solution $w^\eps \in \dot H^1$ to \eqref{stat}, \emph{i.e.} such that 
$$J^\eps (w^\eps )= \min \left\lbrace J^\eps (v) ~:~ v\in \dot H^1, \div v =0 \right\rbrace ,$$ where $J^\eps (v)= \int \left(\frac{1}{2} |\nabla v|^2 - f.v \right) \, dx$. Now, using a penalization method, we change the problem into minimizing 
$$ J^\epsalpha (v) = J^\eps (v) + \frac{1}{2\alpha} \int |\div v|^2 \, dx$$ in the space $\dot H^1$ so that the constraint $\div w^\eps =0$ is integrated to the functional $J^\epsalpha$ to minimize. Let us call the minimizer $w^\epsalpha$. Letting $\alpha $ go to zero in $w^\epsalpha$, we obtain the desired solution $w^\eps$. Since $w^\epsalpha$ minimizes $J^\epsalpha$, we know that it solves
\begin{equation}
 -\Delta w^\epsalpha =  f + \frac{1}{\alpha} \nabla (\div w^\epsalpha) =0 .
\end{equation}
Now, recall that $f \in\dot H^{-1}$ and write
\begin{eqnarray*}
 \int_{\R^2} \left(|\nabla w^\epsalpha|^2 + \frac{1}{\alpha} | \div w^\epsalpha |^2 \right) \, dx& = & \int_{\R^2} f.w^\epsalpha \, dx\\
 & \leq & \|f\|_{\dot H^{-1}} \|\w\|_{\dot H^1} \\
 & \leq & \frac{1}{2} \|f\|_{\dot H^{-1}}^2 + \frac{1}{2} \|\w\|_{\dot H^1}^2 .
\end{eqnarray*}
So we immediately deduce that
\begin{equation} \label{weakcompressibility} \| \div \w \|_{L^2}^2 \leq C \alpha \|f\|_{\dot H^{-1}}^2 = \mathcal{O}(\alpha). \end{equation}
On this basis, we shall consider the equation
 $$ (HNS^{\eps , \alpha}) ~~~ \eps \partial_{tt} u^{\eps , \alpha} + \partial_t u^{\eps , \alpha} - \Delta u^{\eps , \alpha} = -(\u . \nabla)\u + \frac{1}{\alpha} \nabla (\div u^{\eps , \alpha}).$$
 
Due to \eqref{weakcompressibility}, we say that $(HNS^\epsalpha)$ is weakly compressible. Let us point out here that this model reminds of the one studied by Vi\v{s}ik and Fursikov in \cite{vishik} and, later, by Basson \cite{basson} and Lelièvre \cite{lelievre} in order to study statistical solutions to the Navier-Stokes equations. \\

In the next section, we shall prove that $ (HNS^{\eps , \alpha})$ has a finite propagation speed through Picard's fixed point theorem and energy inequalities.

\section{Finite speed of propagation} \label{vitessefinie}
In order to prove that the equation 
\[ (HNS^\epsalpha)~~ \eps \partial_{tt} u^{\eps , \alpha} + \partial_t u^{\eps , \alpha} - \Delta u^{\eps , \alpha} = -(\u . \nabla)\u + \frac{1}{\alpha} \nabla (\div u^{\eps , \alpha})\] 
has a finite propagation speed, let us consider the Helmholtz-Hodge decomposition of $\u$: 
$$\u = \w + \z,$$ 
where $\w = \q \u := \frac{1}{\Delta} \nabla \div \u$ is irrotational and $\z = \p \u := \u - \w$ is divergence-free. We obtain the system
$$
\left\lbrace
\begin{array}{rcl}
 \eps \partial_{tt} \z + \partial_t \z - \Delta \z & = & - \p \left( (\w+\z).\nabla \right) (\w+\z)  \\
 \eps \partial_{tt} \w + \partial_t \w - \frac{\alpha +1}{\alpha}  \Delta \w & = & - \q \left( (\w+\z).\nabla \right) (\w+\z)  
\end{array}
\right.
$$
from which we can deduce a Duhamel's formula for the Cauchy problem
$$ (P) ~~
\left\lbrace
\begin{array}{l}
 (HNS^\epsalpha) \\
 u |_{t=0} = u_0 ~, ~~ \dt u |_{t=0} = u_1
\end{array}
\right. .
$$
In section \ref{existencelocale}, we show local existence through Picard's contraction theorem in a small ball of the complete metric space
\begin{eqnarray*}
X_T & =&  \left\lbrace
	    (u,\partial_t u ) \in \dot H^{\frac{n}{2}+\delta} \cap \dot H^{\frac{n}{2}+\delta-1} (\R^n) \times  \dot H^{\frac{n}{2}+\delta-1} (\R^n)  \right\rbrace .
\end{eqnarray*}
The contractive map argument is detailed in section \ref{existencelocale}. In the following, we shall denote $\u$ by $u$ to alleviate  the notations. Picard's fixed point theorem gives a sequence 
$$ (u^j , \dt u^j ) \stackrel{X_T}{\longrightarrow} (u, \dt u) $$
defined by
\begin{equation} \label{suitepicard}
 \eps \dtt u^{j+1} - \Delta u^{j+1}=- \dt u^j - (u^j.\nabla)u^j + \frac{1}{\alpha} \nabla \div u^j .
\end{equation}
Now, set $\tilde u = \div u$ and apply the $\div$ operator to $(HNS^\epsalpha)$. Doing so, we obtain the following system:
$$ (S) ~~
\left\lbrace
\begin{array}{lcl}
 \eps \dtt \tilde u - \frac{\alpha +1}{\alpha} \Delta \tilde u &=& - \dt \tilde u- \div (u . \nabla) u  \\
 \eps \partial_{tt} u    -   ~~~~~ \Delta u & = & -\partial_t u -(u.\nabla) u + \frac{1}{\alpha} \nabla \tilde u
\end{array}
\right. .
$$
Then the Cauchy problem $(P)$ is equivalent to
$$ (P') ~~
\left\lbrace
\begin{array}{l}
 (S) \\
 u |_{t=0} = u_0 ~,~~ \dt u |_{t=0} = u_1 ~,~~ \tilde u |_{t=0} = \div u_0 ~,~~ \dt \tilde u |_{t=0} = \div u_1
\end{array}
\right. .
$$
\begin{rk}
 As we shall see in the following, this way of writing the equation $(HNS^\epsalpha)$ is convenient for the proof of finite speed of propagation but, unless we assume that $u$ and $\tilde u$ are smooth, we cannot prove directly that the Duhamel's formula related to system $(S)$ is locally contractive (due to the term $\nabla \tilde u$). Besides, we cannot prove that $(HNS^\epsalpha)$ has a finite speed of propagation through the Helmholtz-Hodge decomposition since the operators $\p$ and $\q$ are non-local.
\end{rk}
Moreover, applying $\div$ to \eqref{suitepicard}, we obtain the equation 
$$ \eps \dtt \tilde u^{j+1} -  \frac{\alpha +1}{\alpha}  \Delta \tilde u^{j+1} = - \dt \tilde u^j- \div (u^j . \nabla) u^j .
$$
Then we have the coupled system of equations
$$ (S^j) ~~
\left\lbrace
\begin{array}{lcl}
 \eps \dtt \tilde u^{j+1} -  \frac{\alpha +1}{\alpha}  \Delta \tilde u^{j+1} &=& - \dt \tilde u^j- \div (u^j . \nabla) u^j  \\
 \eps \partial_{tt} u^{j+1}    -  \, ~~~~~~~~~~ \Delta u^{j+1} & = & -\partial_t u^j -(u^j.\nabla) u^j + \frac{1}{\alpha} \nabla \tilde u^j 
\end{array}
\right. .
$$
Now, let us consider the problem
$$ (P^j) ~~
\left\lbrace
\begin{array}{l}
 (S^j) \\
 u^{j+1} |_{t=0} = u_0 \, ~~~~~,~~ \dt u^{j+1} |_{t=0} = u_1 \\
 \tilde u^{j+1} |_{t=0} = \div u_0 ~,~~ \dt \tilde u^{j+1} |_{t=0} = \div u_1 \\
 u^0 ~,~~ v^0 \equiv 0
\end{array}
\right. .
$$
From now on, we assume $u_0$ and $u_1$ to be supported in a ball $B(0,R)$ and we set
$$c_1 = \sqrt{\frac{\alpha +1}{\alpha \eps}} ~,~~ c_2 = \frac{1}{\sqrt{\eps}}.$$ 
Let us recall the standard energies associated to the wave equations in $(S)$:
$$E_{c_1}(\tilde u)(t)=\int_{\Gamma_{c_1 ,t,x}} e_{c_1}(\tilde u)(t,y) \, \d y ~,~~ E_{c_2}(u)(t)=\int_{\Gamma_{c_2 ,t,x}} e_{c_2}(u)(t,y) \, \d y , $$ where the densities are
\begin{eqnarray*}
e_{c_1}(\tilde u)(t,y) &=& \frac{\eps}{2} \, |\partial_t \tilde u (t,y) |^2 + \frac{\alpha+1}{2\alpha}\, |\nabla \tilde u (t,y) |^2 ,\\
e_{c_2}(u)(t,y) &=& \frac{\eps}{2} \, |\partial_t u (t,y) |^2  + ~~~ \frac{1}{2} ~~~ |\nabla u (t,y) |^2 
\end{eqnarray*}
and the cones $\Gamma_{c_i ,t,x}$ are defined by
$$ \Gamma_{c_i ,t,x} = \left\lbrace (s,y) ~:~0\leq s \leq \frac{t}{c_i} ,~ |y-x|\leq t-c_i s \right\rbrace . $$
Notice that, since $c_1>c_2$, we have $\Gamma_{c_1 ,t,x} \subset \Gamma_{c_2 ,t,x}$.
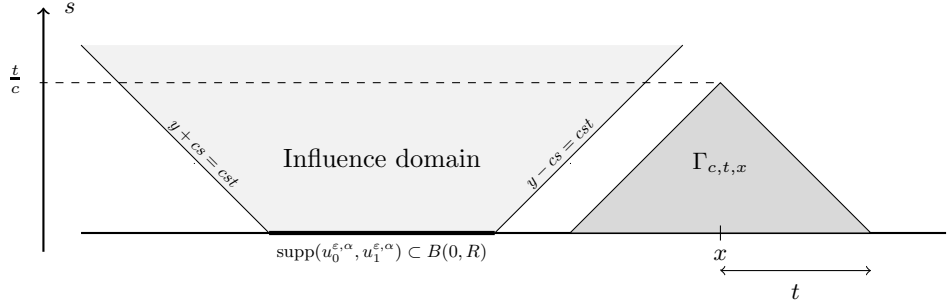
\begin{figure}[ht] \begin{center}
\begin{tikzpicture}[scale=0.5]
	\draw[thick] (-8,0) -- (15,0); 
	\draw[fill=black!5] (-8,5) -- (-3,0) -- (3,0) -- (8,5)  ;
	\draw[<->] (9,-1) -- (13,-1);
	\draw (11,-2) node[anchor=south]{\small{$t$}};
	\draw[ultra thick] (-3,0) -- (3,0);
	\draw (0,-0.95) node[anchor=south, scale=0.75]{\small{$\textrm{supp} (\u_0 , \u_1 ) \subset B(0,R)$}};
	%\draw[fill,red,opacity=0.55] (-0.5,0) -- (1.5,3) -- (3.5, 0) -- cycle; 
	%\draw (1.5,3.95) -- (1.5,3.95) node[anchor=north]{\small{$(x_0 , t_0)$}};
	\draw[fill=black!15] (5,0) -- (13,0) -- (9,4) -- cycle;%(10,3) -- (8,3) -- cycle ;
	\draw (9,0.2) -- (9,-0.2);
	\draw (9,-0.98) -- (9,-0.98) node[anchor=south]{\small{$x$}};
	%\draw[dashed] (8,3) -- (9,4) -- (10,3)  ;
	\draw (9,1.3) node[anchor=south]{\small{$\Gamma_{c,t,x}$}};
	\draw (5,1.85) -- (5,1.85) node[scale=0.7, anchor=south, rotate=45]{\small{$y-c s = cst $}};
	\draw (-5,1.85) -- (-5,1.85) node[scale=0.7, anchor=south, rotate=-45]{\small{$y+c s = cst $}};
	\draw[->,thick] (-9,-0.5) -- (-9,6);
	\draw (-8.7,6) node[anchor=west]{\small{$s$}};
	\draw (-9.3,4.1) node[anchor=east]{\small{$\frac{t}{c}$}};
	\draw[dashed] (-9.1,4) -- (9,4);
	\draw (0,2.5) node[anchor=north]{Influence domain};
\end{tikzpicture}
\caption{The influence and dependence domains for an equation with speed of propagation $c$ and with initial data supported in a ball $B(0,R)$. } \end{center} \end{figure}
Now, assume that we know that $u^j , ~ \dt u^j , ~ \tilde u^j,~ \dt \tilde u^j=0$ on the cone $\Gamma_{c_1 , t,x}$. 
Then, multiplying the first equation in $(S^j)$ by $\dt \tilde u^{j+1}$ and integrating on $\Gamma_{c_1 , t,x}$, we obtain:
\begin{eqnarray*}
 0 & = & \int_{\Gamma_{c_1,t,x}}\!\!\! \left( \eps \dtt \tilde u^{j+1} + \partial_t \tilde u^{j} - \frac{\alpha +1}{\alpha} \Delta \tilde u^{j+1} + \div (u^j .\nabla) u^j \right) \partial_t \tilde u^{j+1} \, \textrm{d}t \, \textrm{d}y \\
 & = & \int_{\Gamma_{c_1,t,x}}\!\!\! \left( \eps \dtt \tilde u^{j+1}  - \frac{\alpha +1}{\alpha} \Delta \tilde u^{j+1} \right) \partial_t \tilde u^{j+1} \, \textrm{d}t \, \textrm{d}y \\
 & = & \int_{\Gamma_{c_1,t,x}} \! \frac{\eps}{2} \dt |\dt \tilde u^{j+1} |^2   - \frac{\alpha +1}{\alpha} \sum_{i=1}^n \partial_{ii} \tilde u^{j+1} . \partial_{t} \tilde u^{j+1} \, \textrm{d}t \, \textrm{d}y  \\
 & = & \int_{\Gamma_{c_1,t,x}} \! \frac{\eps}{2} \dt |\dt \tilde u^{j+1} |^2 - \frac{\alpha +1}{\alpha} \sum_{i=1}^n \left[ \partial_{i} \left( \partial_{i} \tilde u^{j+1} . \partial_{t} \tilde u^{j+1} \right) - \partial_{i} \tilde u^{j+1} . \partial_{ti} \tilde u^{j+1} \right] \, \textrm{d}t \, \textrm{d}y  \\
 & = &  \int_{\Gamma_{c_1,t,x}} \!\!\!\!\!\!\! \dt \left( \frac{\eps}{2} |\dt \tilde u^{j+1} |^2 + \frac{\alpha +1}{2\alpha} | \nabla \tilde u^{j+1} |^2 \right) -\frac{\alpha +1}{\alpha} \sum_{i=1}^n \partial_{i} \left( \partial_{i} \tilde u^{j+1} . \partial_{t} \tilde u^{j+1} \right)  \, \textrm{d}t \, \textrm{d}y .
\end{eqnarray*}
In order to apply the divergence theorem, let us compute the unit outgoing normal $\nu$ to the cone $\Gamma_{c_i,t,x}$. We have
$$ \nu (s,y) = \left\{ \begin{array}{ll}
                 (+1,\mathbf{0}) & \textrm{if} ~ s=t \\
                 (-1,\mathbf{0}) & \textrm{if} ~ s=0 \\
                 \frac{1}{\sqrt{1+c_i^2}} (c_i, \frac{y-x}{|y-x|}) &  \textrm{if}~ 0<s<t.
                \end{array}
\right. 
$$
The divergence theorem yields:
\begin{small}
\[
 E_{c_1}(\tilde u^{j+1})(0) =  E_{c_1}(\tilde u^{j+1}) (t) + \frac{1}{\sqrt{1+c_1^2}} \int_0^{\frac{t}{c_1}} \!\!\! \int_{\mathbb{S}_{t-c_1 s}} \!\!\!\!\!\!\!\!\!\!\!\! c_1 e_{c_1}(\tilde u^{j+1}) (s,y) -\frac{\alpha +1}{\alpha} \sum_{i=1}^n \frac{(y-x)_{i}}{|y-x|} \left( \partial_{i} \tilde u^{j+1} . \partial_{t} \tilde u^{j+1} \right)\!. \] \end{small}
By the Cauchy-Schwarz inequality followed by Young, we obtain:
\[ \int \!\! \sum_{i= 1}^n \frac{(y-x)_{i}}{|y-x|}  \partial_{i} \tilde u^{j+1} . \partial_{t} \tilde u^{j+1} \leq  \sum_{i=1}^n \| \partial_{i} \tilde u^{j+1} \|_{L^2} \| \partial_{t} \tilde u^{j+1} \|_{L^2} \leq  \frac{1}{2 \lambda} \|\nabla \tilde u^{j+1}\|_{L^2}^2 + \frac{\lambda}{2} \|\partial_t \tilde u^{j+1} \|_{L^2}^2 ,\] 
where $\lambda = c_1^{-1}$. Since $\tilde u^{j+1}|_{t=0} = \div u_0$ and $\dt \tilde u^{j+1}|_{t=0} = \div u_1$, we know that \[ \tilde u^{j+1} , \dt \tilde u^{j+1} =0 \textrm{ on } B(x,t) .  \]  So, finally, we have 
$$ E_{c_1}(\tilde u^{j+1}) (t) \leq   E_{c_1}(\tilde u^{j+1})(0) =0 $$ 
and we deduce that $\tilde u^{j+1} ,~ \dt \tilde u^{j+1}=0$ on $\Gamma_{c_1 , t,x}$. \\
\\
Now, in order to handle the second equation in $(S^j)$, we cover the cone $\Gamma_{c_1,x,t}$ with cones of the type $\Gamma_{c_2,t_i,x_i}$ and integrate 
the second equation in $(S^j)$ mutiplied by $\dt u^{j+1}$ on these cones.
\begin{figure}[h] \begin{center}
\begin{tikzpicture}[scale=0.5]
	\draw[thick] (-8,0) -- (8,0); 
	\draw[fill=black!15] (-6,0) -- (6,0) -- (0,6) -- cycle ;
	\draw (0,2.3) node[anchor=south]{\small{$\Gamma_{c_1,t,x}$}};
	\draw[fill=black!5] (-5,0) -- (-3,3) -- (-1,0) -- cycle;
	\draw (-3,0.5) node[anchor=south]{\tiny{$\Gamma_{c_2,t_i,x_i}$}};
	\draw[dashed] (-9.1,3) -- (-3,3) -- (-3,-0.1);
	\draw (-9.2,3) node[anchor=east]{\small{$\frac{t_i}{c_2}$}};
	\draw (-3,-0.98) node[anchor=south]{\small{$x_i$}};
	\draw[fill=black!5] (3,0) -- (4,2) -- (5,0) -- cycle;
	\draw (4,0.5) node[anchor=south]{\tiny{\tiny{$\Gamma_{c_2,t_k,x_k}$}}};
	\draw[dashed] (-9.1,2) -- (4,2) -- (4,-0.1);
	\draw (-9.2,1.9) node[anchor=east]{\small{$\frac{t_k}{c_2}$}};
	\draw (4,-0.98) node[anchor=south]{\small{$x_k$}};
	\draw[->,thick] (-9,-0.5) -- (-9,8);
	\draw (-8.7,8) node[anchor=west]{\small{$s$}};
	\draw (0,0.2) -- (0,-0.2); \draw (0,-0.98) -- (0,-0.98) node[anchor=south]{\small{$x $}};
	%\draw (-9.3,6.1) node[anchor=east]{\small{$t$}};
	%\draw[dashed] (-9.2,6) -- (0.09,6);
\end{tikzpicture}
\end{center}
\caption{Covering the cone $\Gamma_{c_1,x,t}$ with cones of the type $\Gamma_{c_2,t_i,x_i}$.}
\end{figure}
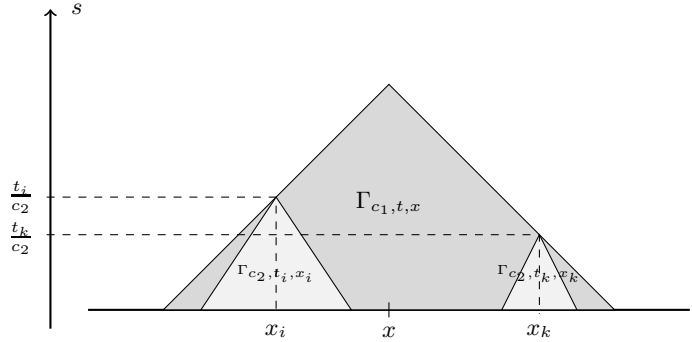
Doing so, we obtain that $u^{j+1}=0$ on $\Gamma_{c_1 , t,x}$ . We have proven that $(HNS^\epsalpha)$ has a finite speed of propagation $c(\epsalpha) \geq c_1 \rightarrow +\infty$ as $\alpha$ goes to $0$ .

\section{Local existence for equation $(HNS^{\epsalpha})$} \label{existencelocale}
% Consider the Cauchy problem 
% $$ \eqref{epsalpha} ~,~~ u|_{t=0} = u_0 \in H^{\frac{3}{2}+\delta}(\mathbb{R}^3) ~,~~ \partial _t u|_{t=0} = u_1  \in H^{\frac{1}{2}+\delta}(\mathbb{R}^3) .$$
% Set $w= \frac{1}{\Delta} \nabla (\div u )$ and $Z=u-w$. Notice that $u=w+Z$ is the Helmholtz decomposition of $u$, where $w$ is the irrotational field and $Z$ the solenoidal one. $w$ is solution to 
% \begin{eqnarray}
% \eps \partial_{tt} w &=& - \partial_t w + \left(1+ \frac{1}{\alpha} \right)\Delta w ,\\ 
% w|_{t=0} &=& w_0:= \frac{1}{\Delta} \nabla (\div u_0 )\in H^{\frac{3}{2}+\delta}(\mathbb{R}^3) ,\nonumber \\
% \partial _t w|_{t=0} &=& w_1:= \frac{1}{\Delta} \nabla (\div u_1 ) \in H^{\frac{1}{2}+\delta}(\mathbb{R}^3) .  \nonumber 
% \end{eqnarray}
% and $Z$ is solution to
% \begin{eqnarray}
% \eps \partial_{tt} Z &=& - \partial_t Z + \Delta Z ,\\ 
% Z|_{t=0} &=& Z_0:= u_0 - w_0 \in H^{\frac{3}{2}+\delta}(\mathbb{R}^3) ,\nonumber \\
% \partial _t w|_{t=0} &=& w_1:= Z_1:= u_1 - w_1  \in H^{\frac{1}{2}+\delta}(\mathbb{R}^3) .  \nonumber 
% \end{eqnarray}
\subsection{Introduction}
Let us consider the Cauchy problem
\begin{equation}  \label{equ}
 \left\lbrace
\begin{array}{l}
 (HNS^\epsalpha) ~~ \eps \partial_{tt} \u + \partial_t \u - \Delta \u - \frac{1}{\alpha} \nabla (\div \u ) = f(\u) \\
 \u |_{t=0}=\u _0 \in H^{\frac{n}{2}+\delta}(\R^n)~,~~ \partial_t \u|_{t=0}=\u_1 \in H^{\frac{n}{2}+\delta-1}(\R^n) 
\end{array}
\right. ,
\end{equation}
where $n=2,3$ and $f(u)= -(u.\nabla )u$. First, let us assume that $f=0$ and split the solution $\u$ to equation $(HNS^\epsalpha)$ into its irrotational part $\w$ and its divergence-free part $\z$: 
$$\u (t,x) = \w (t,x) + \z(t,x) .$$ More precisely, $\w$ and $\z$ are defined as follows:
\begin{eqnarray*}
\w &=& \q \u := \frac{1}{\Delta} \nabla (\div \u) \\ 
\z &=& \p \u := (\mathbf{1}-\q ) \u = \u - \w .
\end{eqnarray*}
Then $\w$ and $\z$ solve the following equations:
\begin{equation}  \label{eqw} 
 \left\lbrace
\begin{array}{l}
 \eps \partial_{tt} \w + \partial_t \w - \left( 1+\frac{1}{\alpha} \right) \Delta \w  = 0\\
 \w|_{t=0}=\w_0:= \q \u_0 ~,~~ \partial_t \w|_{t=0}=\w_1:= \q \u_1  , 
\end{array}
\right.
\end{equation}
 
\begin{equation}  \label{eqz}
 \left\lbrace
\begin{array}{l}
 \eps \partial_{tt} \z + \partial_t \z - \Delta \z  = 0 \\
 \z|_{t=0}=\z_0:= \p \u_0  ~,~~ \partial_t \z|_{t=0}=\z_1:= \p \u_1 .
\end{array}
\right.
\end{equation}
Since \eqref{eqw} is a wave equation, we know that $$\w (t,x) = \aq (t) \w_0 (x) + \bq (t) \w_1 (x) - \int_0^t \bq (t-s) \partial_t \w (s) \, ds,$$ where $\aq$ and $\bq$ are defined as follows.
$$
\aq (t) = \cos \left( \sqrt{\frac{\alpha +1}{\eps \alpha}} t \Lambda \right)~,~~
\bq (t) = \frac{\sin \left( \sqrt{\frac{\alpha +1}{\eps \alpha}} t \Lambda \right)}{\sqrt{\frac{\alpha +1}{\eps \alpha}} \Lambda}.
$$
Similarly, \eqref{eqz} is a wave equation and we have $$\z (t,x) = \ap (t) \z_0 (x) + \bp (t) \z_1 (x) - \int_0^t \bp (t-s) \partial_t \z (s) \, ds,$$ where $\ap$ and $\bp$ are
$$
\ap (t) = \cos \left( \frac{ t \Lambda }{\sqrt{\eps}}  \right)~,~~
\bp (t) = \sqrt{\eps} \frac{\sin \left( \frac{t \Lambda}{\sqrt{\eps }}  \right)}{\Lambda}.
$$
Now, setting $A(t)= \aq (t) \q + \ap (t) \p$ and $B(t) =\bq (t) \q + \bp (t) \p$, we can write Duhamel's formula for the Cauchy problem \eqref{equ} with $f=0$:
$$\phi (\u)(t)= A(t) \u_0 + B(t) \u_1 - \int_0^t B(t-s) \partial_t \u (s) \, ds .$$ 
Adding the source term $f(\u)=-(\u . \nabla) \u$, the formula becomes
$$\phi (\u)(t)= A(t) \u_0 + B(t) \u_1 + \int_0^t B(t-s) \left( f(\u) - \partial_t \u \right)(s) \, ds .$$
\subsection{Contraction argument} \label{picard}
We shall show local existence for \eqref{equ} in the complete metric space 
\begin{eqnarray*}
X_T (a)& =&  \left\lbrace
	    (u,\partial_t u ) \in \dot H^{\frac{n}{2}+\delta} \cap \dot H^{\frac{n}{2}+\delta-1} (\R^n) \times  \dot H^{\frac{n}{2}+\delta-1} (\R^n) ~:~  \right.\\
& & \left. \|u\|_{X_T} := \|u\|_{L^\infty_T \dot H^{\frac{n}{2}+\delta}} + \|u\|_{L^\infty_T \dot H^{\frac{n}{2}+\delta-1}} + \|\partial_t u\|_{L^\infty_T \dot H^{\frac{n}{2}+\delta-1}} \leq a \right\rbrace,
\end{eqnarray*}
with $a>0$ and $0<T<1$ to be chosen later.\\
% A l'aide d'un changement d'\'echelle, on obtient
% \begin{eqnarray*}
%  \aq (t) & = & \cos \left( \sqrt{\frac{\alpha +1}{\alpha \eps}} t \Lambda \right) \\
%  \bq (t) & = & \sqrt{\frac{\eps \alpha}{\alpha +1} } ~ \frac{ \sin \left( \sqrt{\frac{\alpha +1}{\alpha \eps}} t \Lambda \right)}{\Lambda} \\
%  \ap (t) & = & \cos \left( \frac{t}{\sqrt{\eps}} \Lambda \right)\\
%  \bp (t) & = & \sqrt{\eps} ~ \frac{\sin \left( \frac{t}{\sqrt{\eps}} \Lambda \right)}{\Lambda}
% \end{eqnarray*}
Let us estimate $\| \phi (u)  \|_{X_T}$. First, we have
\begin{eqnarray*}
 \| \phi (u) (t) \|_{\dot H^{\frac{n}{2}+\delta}} & \leq & \|\aq \q u_0 \|_{\dot H^{\frac{n}{2}+\delta}} + \| \ap \p u_0 \|_{\dot H^{\frac{n}{2}+\delta}} + \| \bq \q u_1 \|_{\dot H^{\frac{n}{2}+\delta}} + \| \bp \p u_1 \|_{\dot H^{\frac{n}{2}+\delta}}  \\
 & & + \int_0^t \| \bq (t-s) \q \left( f(u) - \partial_t u \right)\|_{\dot H^{\frac{n}{2}+\delta}} \, \d t  \\
 & & + \int_0^t \| \bp (t-s) \p \left( f(u) - \partial_t u \right)\|_{\dot H^{\frac{n}{2}+\delta}} \, \d t \\
 & \leq & 2 \|u_0\|_{\dot H^{\frac{n}{2}+\delta}} + \left(1 + \sqrt{\frac{\alpha}{\alpha+1}} \right) \sqrt{\eps} \|u_1\|_{\dot H^{\frac{n}{2}+\delta-1}} \\
 & & + \sqrt{\frac{\eps \alpha}{\alpha +1}} \int_0^t \left\Vert ~ \frac{\sin \left( \sqrt{\frac{\alpha +1}{\alpha \eps}} (t-s) \Lambda \right)}{\Lambda} ~ \q \left[ (u.\nabla) u - \partial_t u \right] (s)  \right\Vert_{\dot H^{\frac{n}{2}+\delta}}  \d s  \\
% & & + \sqrt{\frac{\eps \alpha}{\alpha +1}} \int_0^t \left\Vert ~ \frac{\sin \left( \sqrt{\frac{\alpha +1}{\alpha \eps}} (t-s) \Lambda \right)}{\Lambda} ~ \q \partial_t u(s) ~ \right\Vert_{\dot H^{\frac{n}{2}+\delta}} \d s + \\
 & & + \sqrt{\eps} \int_0^t \left\Vert ~ \frac{\sin \left( \frac{t-s}{\sqrt{\eps}} \Lambda \right)}{\Lambda} ~ \p  \left[ (u.\nabla) u - \partial_t u \right](s) \right\Vert_{\dot H^{\frac{n}{2}+\delta}} \d s + \\
 %& & + \sqrt{\eps} \int_0^t \left\Vert ~ \frac{\sin \left( \frac{t-s}{\eps} \Lambda \right)}{\Lambda} ~ \p \partial_t u(s) ~ \right\Vert_{\dot H^{\frac{n}{2}+\delta}} \d s + \\
 & = & 2 \|u_0\|_{\dot H^{\frac{n}{2}+\delta}} + \left(1 + \sqrt{\frac{\alpha}{\alpha+1}} \right) \sqrt{\eps} \|u_1\|_{\dot H^{\frac{n}{2}+\delta-1}}  + \MakeUppercase{\romannumeral 1} + \MakeUppercase{\romannumeral 2}.
\end{eqnarray*}
In order to estimate $\MakeUppercase{\romannumeral 1}$, recall that $\div u \neq 0$, so that we have
\[(u.\nabla) u = \sum_{i=1}^n \partial_i \left( u_i u \right) - \div u ~ u .\]
Using that $\dot H^{\frac{n}{2}+\delta} \cap \dot H^{\frac{n}{2}+\delta -1} (\R^n)$ is an algebra, we obtain the inequality
\[\|\partial_i \left( u_i u \right) \|_{\dot H^{\frac{n}{2}+\delta -1}} \leq \| u_i u  \|_{\dot H^{\frac{n}{2}+\delta }} \leq 2C \|u\|_{\dot H^{\frac{n}{2}+\delta }} \|u\|_{\dot H^{\frac{n}{2}+\delta -1}} .\]
Now, the following Lemma allows to estimate the remaining term.
\begin{lemme} \label{div_lp}
 Let $u \in \dot H^{\frac{n}{2}+\delta} \cap \dot H^{\frac{n}{2}+\delta -1} (\R^n)$ and $\div u \neq 0$. Then we have
 \[ \| \div u ~ u \|_{\dot H^{\frac{n}{2}+\delta -1}} \leq 2C  \|u\|_{\dot H^{\frac{n}{2}+\delta}} \|u\|_{L^\infty} .  \]
\end{lemme}
\proof Applying the paraproduct decomposition (see Lemma \ref{paraproduct} in the appendix), we can write
\begin{eqnarray*}
 \| \partial_i u ~ u \|_{\dot H^{\frac{n}{2}+\delta -1}} & \leq & \sum_p 2^{p\left(\frac{n}{2}+\delta -1 \right)} \| \Delta_p \partial_i u S_{p+1} u \|_{L^2} + \sum_q 2^{q\left(\frac{n}{2}+\delta -1 \right)} \| \Delta_q u S_{p+1} \partial_i u \|_{L^2} \\
 & \leq & \sum_p 2^{p\left(\frac{n}{2}+\delta -1 \right)} \| \Delta_p \partial_i u \|_{L^2} \|S_{p+1} u \|_{L^\infty} + \sum_q 2^{q\left(\frac{n}{2}+\delta -1 \right)} \| \Delta_q u \|_{L^2} \|S_{p+1} \partial_i u \|_{L^\infty} \\
 & \leq & \sum_p 2^{p\left(\frac{n}{2}+\delta  \right)}  \| \Delta_p u \|_{L^2} \| u \|_{L^\infty} + \sum_q 2^{q\left(\frac{n}{2}+\delta  \right)} \| \Delta_q u \|_{L^2} \| u \|_{L^\infty} \\
 & \leq & 2C\|u\|_{\dot H^{\frac{n}{2}+\delta}} \|u\|_{L^\infty}.
\end{eqnarray*}
The Lemma is thereby proved. \hfill $\Box$ \\
\\
Now, using the interpolation estimate $\|u\|_{L^\infty} \leq C \|u\|_{\dot H^{\frac{n}{2}+\delta -1}}^\delta \|u\|_{\dot H^{\frac{n}{2}+\delta }}^{1-\delta} $ followed by a Young inequality, we obtain that
\[ \| \div u ~ u \|_{\dot H^{\frac{n}{2}+\delta -1}}  \leq 2C \|u\|_{\dot H^{\frac{n}{2}+\delta -1}} \|u\|_{\dot H^{\frac{n}{2}+\delta}} +\|u\|_{\dot H^{\frac{n}{2}+\delta}}^2 . \]
Finally, we have
\[ \| (u.\nabla) u \|_{\dot H^{\frac{n}{2}+\delta -1}} \leq 4C \| u \|_{\dot H^{\frac{n}{2}+\delta } \cap \dot H^{\frac{n}{2}+\delta -1}}^2 .  \]
and the integral $\MakeUppercase{\romannumeral 1}$ estimates
\begin{eqnarray*}
 \MakeUppercase{\romannumeral 1} & \leq & \sqrt{\frac{\eps \alpha}{\alpha +1}} \int_0^t \left( \| \q (u.\nabla) u (s)\|_{\dot H^{\frac{n}{2}+\delta -1}} + \| \q \partial_t u (s) \|_{\dot H^{\frac{n}{2}+\delta -1}}  \right) \, \d s \\
 & \leq &  \sqrt{\frac{\eps \alpha}{\alpha +1}} \int_0^t \left( \| (u.\nabla) u (s)\|_{\dot H^{\frac{n}{2}+\delta -1}} + \| \partial_t u (s) \|_{\dot H^{\frac{n}{2}+\delta -1}}  \right) \, \d s \\
 & \leq & \sqrt{\frac{\eps \alpha}{\alpha +1}}  \int_0^t \left( 4C \| u (s)\|_{\dot H^{\frac{n}{2}+\delta } \cap \dot H^{\frac{n}{2}+\delta -1}}^2 + \| \partial_t u (s) \|_{\dot H^{\frac{n}{2}+\delta -1}}  \right) \, \d s \\
 & \leq & \sqrt{\frac{\eps \alpha}{\alpha +1}}  \left( 4C T \| u (s)\|_{L^\infty_T \left( \dot H^{\frac{n}{2}+\delta } \cap \dot H^{\frac{n}{2}+\delta -1} \right)}^2 + T \| \partial_t u (s) \|_{L^\infty_T \dot H^{\frac{n}{2}+\delta -1}}  \right) \\
 & \leq & \sqrt{\frac{\eps \alpha}{\alpha +1}} T ( 4C a^2 + a  ) .\\
\end{eqnarray*}
Similarly, we obtain 
$$ \MakeUppercase{\romannumeral 2} \leq \sqrt{\eps} T ( 4C a^2 + a).$$
So, finally, we have the inequality
\begin{eqnarray*}
 \| \phi (u) \|_{L^\infty_T \dot H^{\frac{n}{2}+\delta}} & \leq & 2 \|u_0\|_{\dot H^{\frac{n}{2}+\delta}} + \left(1 + \sqrt{\frac{\alpha}{\alpha+1}} \right) \sqrt{\eps} \|u_1\|_{\dot H^{\frac{n}{2}+\delta-1}} + \\
 & & +\sqrt{\eps} \left( 1+ \sqrt{\frac{\alpha}{\alpha +1}} \right) aT (1 + 4C a).
\end{eqnarray*}
Analogous estimates using that $\left\vert \frac{\sin x}{x} \right\vert \leq 1$ give
$$\| \phi (u) \|_{L^\infty_T \dot H^{\frac{n}{2}+\delta-1}} \leq  2 \|u_0\|_{\dot H^{\frac{n}{2}+\delta -1}} + 2T \|u_1\|_{\dot H^{\frac{n}{2}+\delta -1}} +  2aT^2 (1 + 4C a). $$
Then, let us compute the time derivative of $\phi (v)$.
\begin{eqnarray*}
 \partial_t \phi (v)(t) & = &  - \sqrt{\frac{\alpha +1}{\alpha \eps}} \Lambda \sin \left( \sqrt{\frac{\alpha +1}{\alpha \eps}} t \Lambda \right) \q v_0 - \frac{\Lambda}{\sqrt{\eps}} \sin \left( \frac{t}{\sqrt{\eps}} \Lambda \right) \p v_0 +\\
 & & + \cos \left( \sqrt{\frac{\alpha +1}{\alpha \eps}}t \Lambda \right) \q v_1 \hspace*{1.7cm}+ \cos \left( \frac{t}{\sqrt{\eps}} \Lambda \right) \p v_1 + \\
 & & + \int_0^t \cos \left( \sqrt{\frac{\alpha +1}{\alpha \eps}} (t-s) \Lambda \right) \q \left( v.\nabla v - \partial_t v \right)(s) \, \d s + \\
 & & + \int_0^t \cos \left( \frac{t-s}{\sqrt{\eps}} \Lambda \right) \p \left( v.\nabla v - \partial_t v \right)(s) \, \d s .
\end{eqnarray*}
We have 
$$ \| \partial_t \phi (u) \|_{L^\infty_T \dot H^{\frac{n}{2}+\delta-1}} \leq \frac{1}{\sqrt{\eps}} \left( 1+ \sqrt{\frac{\alpha +1}{\alpha}} \right) \|u_0\|_{\dot H^{\frac{n}{2}+\delta}} + 2 \|u_1\|_{\dot H^{\frac{n}{2}+\delta-1}} + 2aT (1+4Ca).$$
Finally, we obtain
\begin{eqnarray*}
  \| \phi (u) \|_{X_T} & \leq & \left( 2+ \frac{1}{\sqrt{\eps}} + \sqrt{\frac{\alpha +1}{\alpha \eps}} \right) \|u_0\|_{\dot H^{\frac{n}{2}+\delta}} + 2 \|u_0\|_{\dot H^{\frac{n}{2}+\delta-1}} + \\
  & & + \left( 2 + \sqrt{\eps} + \sqrt{\frac{\alpha \eps }{\alpha +1}} \right) \|u_1\|_{\dot H^{\frac{n}{2}+\delta-1}} + 2T \|u_1\|_{\dot H^{\frac{n}{2}+\delta-1}} +\\
  & & + \left( 2 + 2T + \sqrt{\eps} + \sqrt{\frac{\alpha \eps }{\alpha +1}} \right) aT (1+4Ca).
\end{eqnarray*}
Let us set $$ \frac{a}{2} =\left( 2+ \frac{1}{\sqrt{\eps}} + \sqrt{\frac{\alpha +1}{\alpha \eps}} \right) \|u_0\|_{\dot H^{\frac{n}{2}+\delta}} + 2 \|u_0\|_{\dot H^{\frac{n}{2}+\delta-1}} + \left( 2 + \sqrt{\eps} + \sqrt{\frac{\alpha \eps }{\alpha +1}} \right) \|u_1\|_{\dot H^{\frac{n}{2}+\delta-1}}.
$$
So we have
$$ \| \phi (u) \|_{X_T} \leq \frac{a}{2} + \frac{a T}{2 + \sqrt{\eps} + \sqrt{\frac{\alpha \eps }{\alpha +1}}} + \left( 2 + 2T + \sqrt{\eps} + \sqrt{\frac{\alpha \eps }{\alpha +1}} \right) aT (1+4Ca).
$$
Finally, we have the following bound on the local existence time:
$$T \leq \frac{C}{1  +  \left[ \left( C_\eps + \sqrt{\frac{\alpha +1 }{\alpha \eps}} \right)  \|\u_0\|_{\dot H^{\frac{n}{2}+\delta}} + 2 \|\u_0\|_{\dot H^{\frac{n}{2}+\delta -1}} + \left( \tilde C_\eps + \sqrt{\frac{\alpha \eps}{\alpha +1}} \right) \|\u_1\|_{\dot H^{\frac{n}{2}+\delta -1}} \right]} \cdot$$
In order to prove that the solutions obtained in this section are global, we shall prove that the denominator 
$$ \left( C_\eps + \sqrt{\frac{\alpha +1 }{\alpha \eps}} \right)  \|\u_0\|_{\dot H^{\frac{n}{2}+\delta}} + 2 \|\u_0\|_{\dot H^{\frac{n}{2}+\delta -1}} + \left( \tilde C_\eps + \sqrt{\frac{\alpha \eps}{\alpha +1}} \right) \|\u_1\|_{\dot H^{\frac{n}{2}+\delta -1}} $$ remains bounded for all fixed $\alpha$ and $\eps$. \\
\\
In subsections \ref{glob2d} and \ref{glob3d},we will globalize the solutions obtained in this section and then, in subsections \ref{conv2d} and \ref{conv3d}, we will prove that they converge towards the solutions to $(HNS^\eps)$ as $\alpha$ goes to $0$.

\section{The 2D case} \label{2d}

\subsection{Preliminary estimates} \label{rappels2d}
\newcommand{\be}{\begin{equation}}
\newcommand{\ee}{\end{equation}}                  
In this part, we shall recall the regularity results we have on the solution $u^\eps$ to $(HNS^\eps)$. %All the proofs are in \cite{article}.
\begin{lemme} Let $T>0$ and $u^\eps \in L^\infty_T (\dot H^{1+\delta} \cap \dot H^{\delta})(\R^2)^2$ be the global solution to $(HNS^\eps)$ with initial data $(u_0^\eps , u^\eps_1 ) \in H^{1+\delta} \times H^{\delta} (\R^2)^2$. Then we have
 $$u^\eps \in L^\infty_T \dot H^{1+\delta} \cap L^2_T \dot H^{1+\delta} \cap L^2_T L^2 \cap L^\infty_T L^2 \cap L^2_T \dot H^1$$
 and
 $$\dt u^\eps \in L^2_T L^2.$$
\end{lemme}
We sketch here the proof of this lemma since all the details can be found in \cite{article}. In the following, we shall denote by $C$ all the constants, even those depending on $T$.
\proof First, let us introduce the energy
$$E^\delta_\eps (t):= \frac{1}{2} \| u^\eps + \eps \dt u^\eps \|_{\dot H^\delta}^2 + \frac{\eps ^2}{2} \| \dt u^\eps \|_{\dot H^\delta}^2 + \eps \| \nabla u^\eps \|_{\dot H^\delta}^2 ,$$
for non-negative $\delta$. Then, according to \cite{article}, we know that
$$ \exists C>0 ~:~ \forall t \geq 0 ~,~~ E_\eps^\delta (t) \leq C \eps^{-\delta} . $$
From this inequality, we immediately deduce that
\begin{equation} \eps \| u^\eps \|_{L^\infty_T \dot H^{1+\delta} }  + \eps \| u^\eps \|_{L^2_T \dot H^{1+\delta}} \leq C \eps^{-\delta} ~ \textrm{and} ~ \|u^\eps\|_{L^\infty_T L^\infty} = o\left( \frac{1}{\sqrt{\eps}} \right).   \end{equation}
Now, let us compute the time derivative of $ E^0_\eps $. We have 
$$
 \frac{d}{dt}E^0_\eps (t) + \eps \int_{\mathbb{R}^2}|\partial _t u^\eps + \nabla : (u^\eps \otimes u^\eps )|^2 \,dx + \int_{\mathbb{R}^2} \left( |\nabla u^\eps |^2 - \eps |\nabla (u^\eps \otimes u^\eps)|^2 \right)\,dx =0.
$$
Due to the control of the norm $\|u^\eps \|_{L^\infty_T L^\infty}$ by $\dfrac{1}{C\sqrt{\eps}}$ for any $C$ provided that $\eps$ is small enough, the last term in the left hand side is lower bounded by $\int_{\mathbb{R}^2} \frac{1}{2} |\nabla u^\eps |^2 \,\d x$. Now, integrating in time, we have
$$
 E^0_\eps (T) + \int_0^T \!\!\! \int_{\mathbb{R}^2}\eps |\partial _t u^\eps + \nabla : (u^\eps \otimes u^\eps )|^2 \,\d x\,\d t + \int_0^T \!\!\!\int_{\mathbb{R}^2} \frac{1}{2} |\nabla u^\eps |^2 \,\d x \, \d t \leq E^0_\eps (0) \leq C .
$$
So, we have that $\eps \| \partial _t u^\eps + \nabla : (u^\eps \otimes u^\eps )\|_{L^2_T L^2}^2  \leq C$ and 
\begin{equation} \frac{1}{2} \|u^\eps \|_{L^2_T \dot H^1}^2  \leq C . \end{equation}
The later yields
$$
 \eps \| \nabla : (u^\eps \otimes u^\eps )\|_{L^2_T L^2}^2  \leq C \eps \|u^\eps\|_{L^\infty_T L^\infty}^2 \| u^\eps \|_{L^2_T \dot H^1}^2  \leq C .
$$
Consequently, we have
\begin{equation}
\| \sqrt{\eps} \partial_t u^\eps \|_{L^2_T  L^2}^2  \leq  2 \eps \|\partial _t u^\eps + \nabla : (u^\eps \otimes u^\eps )\|_{L^2_T L^2}^2 + 2 \eps \| \nabla : (u^\eps \otimes u^\eps )\|_{L^2_T L^2}^2  \leq  C.
\end{equation}
We have thereby proven that 
$$\sqrt{\eps} \partial_t u^\eps \in L^2_T  L^2 \textrm{ uniformly in } \eps.$$ 
Finally, notice that $\|u^\eps\|_{L^2}^2 \leq 2 E_\eps^0 (t) \leq 2 C_0$, so that \begin{equation} u^\eps \in L^2_T L^2 .\end{equation}
\subsection{Globalization} \label{glob2d}
As in \cite{article}, let us define the energy 
$$E_\epsalpha^\delta (t) = \int_{\R^2} \frac{1}{2}|\Lambda^\delta (\u + \eps \dt \u)|^2 + \frac{\eps ^2}{2} |\Lambda^\delta \dt \u |^2 + \eps |\Lambda^\delta \nabla \u |^2 + \frac{\eps}{\alpha} |\Lambda^\delta \div  \u|^2 .$$
In view of the dependence of the local time existence on the initial data, we know that proving that $E_\epsalpha^\delta$ is bounded yields the global existence for $(HNS^\epsalpha)$. We shall first show that it is true on a time interval $[0,T)$ then prove that $T=+\infty$.\\
\\
First, let us point out that $\dot{H}^{1+\delta} \cap L^\infty (\R^2)$ is an algebra and that the product estimate
\begin{equation}
\label{estimation_douce}
 \left\Vert fg \right\Vert_{\dot{H}^{1+\delta}(\R^2)} \leq C_1 \left( \left\Vert f \right\Vert_{\dot{H}^{1+\delta}} \left\Vert g \right\Vert_{\infty} + \left\Vert g \right\Vert_{\dot{H}^{1+\delta}} \left\Vert f \right\Vert_{\infty} \right).
\end{equation}
holds (see Proposition \ref{tame} in the appendix or \cite{alinhac}) for all functions $f,g \in\dot{H}^{1+\delta} \cap L^\infty (\R^2)$. Moreover, we know that the homogeneous Besov\footnote{For definitions and properties of the Besov spaces, see the book by P.-G. Lemarié-Rieusset \cite{pglr}} space $\dot{B}_{2,1}^1(\mathbb{R}^2)$ embeds into $L^\infty (\mathbb{R}^2)$ and, interpolating, we obtain 
\begin{equation}
\label{interp_besov}
 \|f\|_{\infty} ~ \leq ~ \tilde C \|f\|_{\dot{B}_{2,1}^1} ~ \leq ~ C_2 \left\Vert f \right\Vert_{\dot{H}^\delta}^\delta . \left\Vert f \right\Vert_{\dot{H}^{1+\delta}}^{1-\delta}.
\end{equation}
Finally, notice that $\div \u \neq 0$, so that we have
 \be \label{identitediv} (\u.\nabla )\u = \sum_{i=1}^2 \u_i \partial _i \u  = \sum_{i=1}^2 \partial _i (\u_i . \u) - \div \u \times \u. \ee
 but we still can prove the estimate
 \begin{equation} \label{nonlineariteH} \| (\u.\nabla )\u \|_{\dot H^\delta} \leq C_3 \|\u\|_{L^\infty} \|\u\|_{\dot H^{1+\delta}} \end{equation}
 using that $\dot H^{1+\delta}\cap L^\infty (\R^2)$ is an algebra for the first term in \eqref{identitediv} and due to Lemma \ref{div_lp} for the second term.\\ %Indeed, the first term in the right hand side of \eqref{identitediv} easily estimates due to proposition \ref{tame} (which applies for $\delta >0$) and we have
%\begin{small}
% \begin{eqnarray*}
%  \| \partial_i \u ~ \u \|_{\dot H^\delta} & \leq & \sum_{p\geq -1} 2^{p \delta} \| \Delta_p \partial_i \u ~ S_{p+1} \u \|_{L^2} + \sum_{q \geq -1} 2^{q \delta}\| \Delta_q \u ~ S_{q} \partial_i \u \|_{L^2} \\
%  & \leq & \sum_{p\geq -1} 2^{p \delta} \| \Delta_p \partial_i \u \|_{L^2} \| S_{p+1} \u \|_{L^\infty} + \sum_{q \geq -1} 2^{q \delta} \| \Delta_q \u \|_{L^2} \|S_{q} \partial_i \u \|_{L^\infty} \\
%  & \leq & \sum_{p\geq -1}   c_p \|\u \|_{\dot H^{1+\delta}} \|  \u \|_{L^\infty} + \sum_{q \geq -1} 2^{-q} \tilde c_q \| \u \|_{\dot H^{1+\delta}} 2^q \|S_{q}  \u \|_{L^\infty} \\
%  & \leq &  C  \|\u \|_{\dot H^{1+\delta}} \|  \u \|_{L^\infty} ,
% \end{eqnarray*}
% \end{small}
%where $(c_p) , (\tilde c_q) \in \ell^2$. 

\newcommand{\tmax}{T^{\textrm{max}}}
Now, let us define $\tmax$ the maximal existence time of $(HNS^\epsalpha)$ and prove the following lemma.
\begin{lemme} \label{lemme_controle_energie}
Assume the following, when $\eps$ goes to zero :
$$
 (H) \left\lbrace \begin{array}{ll}
                   i)  &  \eps^{\frac{1+\delta}{2}} \|\u_0\|_{\dot H^{1+\delta}} + \eps^{\frac{\delta}{2}} \|\u_0\|_{\dot H^{\delta}} = o(1) \\
		   ii) &  \eps^{\frac{1}{2}} \|\u_0\|_{\dot H^{1}} + \eps  \|\u_1\|_{L^2} =o(1) .
                  \end{array}
\right.
$$
Let us define $0 \leq T \leq \tmax$ by
\be
 T = \sup \left\lbrace 0 \leq \tau \leq \tmax ~ :~ \forall ~ t \in [0, \tau) , ~\|\u  (t)\|_{L^\infty} < \dfrac{1}{2 C_3 \sqrt{\eps}} \right\rbrace .
\ee
Then, for $\eps$ small enough, there exists a large number $N$, depending only on $\delta$ and $\|\u_0\|_{L^2}$ (which is arbitrary), such that, for all $0 \leq t < T$,
\begin{equation} \label{controle_energie}
 E^\delta_\epsalpha (t) \leq E^\delta_\epsalpha (0) \left(2 \|\u_0\|_{L^2}^2 +1\right)^N .
\end{equation}
\end{lemme}
\proof
Let us compute the time derivative of $E_\epsalpha^\delta$:
\begin{eqnarray*}
 \frac{\textrm{d}}{\textrm{d}t}E_\epsalpha^\delta (t) & = & \eps \|(\u.\nabla) \u\|_{\dot H^\delta}^2 - \| \nabla \u\|_{\dot H^\delta}^2  -\int_{\R^2} \Lambda^\delta \u . \Lambda^\delta (\u.\nabla ) \u \, \textrm{d}x   - \\
 & & - \eps \| \dt \u + (\u.\nabla )\u\|_{\dot H^\delta}^2 - \frac{1}{\alpha} \| \div  \u\|_{\dot H^\delta}^2  \\
 & \leq & \eps \|(\u.\nabla) \u\|_{\dot H^\delta}^2 - \| \nabla \u\|_{\dot H^\delta}^2 - \int_{\R^2} \Lambda^\delta \u . \Lambda^\delta (\u.\nabla ) \u \, \textrm{d}x .
 \end{eqnarray*}
 Using \eqref{nonlineariteH}, the derivative of the energy estimates as follows:
\begin{eqnarray*}
  \frac{\textrm{d}}{\textrm{d}t}E_\epsalpha^\delta (t) & \leq &  \left( C_3^2 \eps \|\u\|_{\infty}^2 - 1 \right) \| \u\|_{\dot H^{1+\delta}}^2  -\int_{\R^2} \Lambda^\delta \u . \Lambda^\delta (\u.\nabla ) \u \, \textrm{d}x .
\end{eqnarray*}
Since we assume $i)$ and using inequality \eqref{interp_besov}, we can write, for $\eps$ small enough,
\begin{equation}
\label{controle_norme_infinie}
 \| u_0^\epsalpha \|_{L^\infty} \leq \dfrac{1}{2 C_3 \sqrt{\varepsilon}}  .
\end{equation}
Now, by continuity of the (local) solution $u^\epsalpha$ with respect to $t$, we deduce that $T>0$ and that the inequality 
\begin{eqnarray*}
  \frac{\textrm{d}}{\textrm{d}t}E_\epsalpha^\delta (t) & \leq & - \frac{1}{4} \| \u\|_{\dot H^{1+\delta}}^2 + C \| \u \|_{\dot H^\delta} \| \u \|_{\infty} \| \u \|_{\dot H^{1+\delta}}
\end{eqnarray*}
holds on $[0,T)$. Then, using the interpolation inequalities
\be
 \|\u \|_{\dot{H}^\delta} \leq C_4 \|\u \|_{L^2}^{1-\delta} \lVert \u \rVert_{\dot{H}^1}^\delta 
\ee
and \eqref{interp_besov}, we obtain the estimate
\[
 \frac{d}{dt} E^\delta_\epsalpha (t) ~ \leq ~ -\dfrac{1}{4} \|\u \|_{\dot{H}^{1+\delta}}^2 + C  \|\u \|_2^{1-\delta} \left( \lVert \u  \rVert_{\dot{H}^1} \lVert \u  \rVert_{\dot{H}^\delta} \right)^\delta \lVert \u  \rVert_{\dot{H}^{1+\delta}}^{2-\delta}.
\]
Finally, a Young inequality yields
\begin{equation} \label{controle_delta}
 \frac{\textrm{d}}{\textrm{d}t}E_\epsalpha^\delta (t) \leq C \|\u\|_2^{2 \frac{1-\delta}{\delta}} \|\u\|_{\dot H^1}^2 E_\epsalpha^\delta(t).
\end{equation}
In order to show that $E_\epsalpha^\delta $ is bounded, we shall use the decay of $E_\epsalpha^0$. So, let us estimate its time derivative:
$$ \frac{\textrm{d}}{\textrm{d}t}E_\epsalpha^0 (t) \leq  \eps \|(\u.\nabla) \u\|_{L^2}^2 - \| \u\|_{\dot H^1}^2  -\int_{\R^2}  \u .  (\u.\nabla ) \u \, \textrm{d}x - \frac{1}{\alpha} \| \div \u \|_{L^2}^2   . $$
First, by Lemma \ref{div_lp}, we immediately have
$$ \| (\u.\nabla )\u \|_{L^2} \leq C_3 \|\u\|_{L^\infty} \|\u\|_{\dot H^1}.$$ 
Besides, by integrations by parts, we obtain
$$\int_{\R^2} \u.(\u.\nabla )\u \, \textrm{d}x  = - \frac{1}{2} \int_{\R^2} \u . \, \div \u ~ \u = - \frac{1}{2} \int_{\R^2} \div \u ~ |\u|^2.$$
Recall that $\div \u \in L^2(\R^2)$ and $\u \in L^2 \cap \dot H^1 \subset \dot H^{\frac{1}{2}} \subset L^4 (\R^2)$. Thus the integral estimates
\begin{eqnarray*}
\lvert \int_{\R^2}\u.(\u.\nabla )\u \, \textrm{d}x \rvert & = & \lvert \frac{1}{2} \int_{\R^2} \div \u ~ |\u|^2 \rvert \\
& \leq &  \frac{1}{2} \| \div \u \|_{L^2}  \||\u|^2 \|_{L^2} \\
& \leq & \frac{1}{2} \| \div \u \|_{L^2}  \|\u \|_{L^4}^2 \\
%& \leq & \tilde C \| \div \u \|_{L^2} \|\u \|_{\dot H^{\frac{1}{2}}}^2 \\
& \leq & \frac{K}{2} \| \div \u \|_{L^2} \|\u \|_{L^2}  \|\u \|_{\dot H^1} \\
& \leq & \frac{K^2}{8} \| \div \u \|_{L^2}^2 \|\u \|_{L^2}^2 + \frac{1}{2} \|\u \|_{\dot H^1}^2,
\end{eqnarray*}
where we have used the Gagliardo-Nirenberg inequality
\[
 \|f\|_{L^4}^2 \leq K \|f\|_{L^2} \|f\|_{\dot H^1}
\]
for $f\in H^1(\R^2)$.
We therefore obtain the following estimate on $E_\epsalpha^0$ on $[0,T)$ for $\eps$ small enough:
\begin{eqnarray*}
 \frac{\textrm{d}}{\textrm{d}t}E_\epsalpha^0 (t) &\leq & \eps \|(\u.\nabla) \u\|_{L^2}^2 - \| \u\|_{\dot H^1}^2  -\int_{\R^2}  \u . (\u.\nabla ) \u \, \textrm{d}x  - \frac{1}{\alpha} \| \div \u\|_{L^2}^2\\
 %& \leq & (2 \eps \|\u\|_{L^\infty}^2 -1) \|\u\|_{\dot H^1}^2 - \frac{1}{\alpha} \| \div \u \|_{L^2}^2 + \frac{C^2}{2} \| \div \u \|_{L^2}^2 \|\u \|_{L^2}^2 + \frac{1}{2} \|\u \|_{\dot H^1}^2 \\
 & \leq & \left(C_3^2 \eps \|\u\|_{L^\infty}^2 -\frac{1}{2}\right) \|\u\|_{\dot H^1}^2 + \left(\frac{K^2}{8}   \|\u \|_{L^2}^2 -\frac{1}{\alpha} \right)  \| \div \u \|_{L^2}^2 \\
 & \leq & -\frac{1}{4} \|\u\|_{\dot H^1}^2 + \left(\frac{K^2}{8}   \|\u \|_{L^2}^2 -\frac{1}{\alpha} \right)  \| \div \u \|_{L^2}^2 .
\end{eqnarray*}
Moreover, assuming that $ \alpha \leq \dfrac{2}{K^2 \|\u_0\|_{L^2}^2 }$, we have
$$ \frac{\textrm{d}}{\textrm{d}t}E_\epsalpha^0 (0) \leq -\frac{1}{4} \|\u_0\|_{\dot H^1}^2 .$$
Besides, the smallness assumptions $ii)$ on the initial data yield
$$
\|\u(t)\|_{L^2}^2 ~ \leq ~ 2 E_\epsalpha^0 (t) ~ \leq ~ 2 E_\epsalpha^0 (0) ~ \stackrel{ii)}{\leq}~ 4 \|\u_0\|_{L^2}^2 
$$ 
and $E_\epsalpha^0$ is therefore decreasing on $[0, T)$. \\
Let us now define the functional $\mathcal{E}_\epsalpha^\delta := E_\epsalpha^\delta (1+E_\epsalpha^0)^N$.  One can easily check that
\[
 \frac{\d}{\d t} \mathcal{E}_\epsalpha^\delta \leq E_\epsalpha^\delta (1+E_\epsalpha^0)^{N-1} \| \u \|_{\dot H^1}^2 \left( C \|\u_0\|_{L^2}^{2\frac{1-\delta}{\delta}} (1+2 \| \u_0 \|_{L^2}^2) - \frac{N}{4} \right).
\]
So, if the positive integer $N=N(\delta, \|\u_0\|_2)$ is large enough, $\mathcal{E}_\epsalpha^\delta$ decreases and $E_\epsalpha^\delta (t)$ satisfies the estimate \eqref{controle_energie} for all $t \in [0,T)$.  \qed \\

The aim of the following lemma is to ensure the control of $\|\u (t)\|_\infty$ throughout the time, so that we can reiterate the reasoning.
\begin{lemme} \label{lemme2}
Assume that
\be \label{hyp}
 \eps^{1+\frac{\delta}{2}} \|\u_1\|_{\dot H^\delta} \longrightarrow 0~, ~~\eps \rightarrow 0
\ee
in addition to the assumptions $(H)$ in Lemma \ref{lemme_controle_energie}. Then the inequality
\begin{equation}
\label{controle_initial}
 E^\delta_\epsalpha (0) ~ < ~ C_\delta   (2 \| u_0^\epsalpha \|_2^2 +1)^{-N}.
\end{equation}
holds for $\eps$ small enough. Moreover, for all $t \in [0,T)$,
\begin{equation}  \label{ongardelecontrole}
\|u^\epsalpha \|_{L^\infty} ~\leq ~ \frac{1}{4 C_3 \sqrt{\eps}} .
\end{equation}
\end{lemme}
We skip the proof of this lemma since it can be found in \cite{article}.
\begin{rk}
 Notice that, under the assumptions \eqref{size} in Theorem \ref{th_alpha}, the conditions $(H)$ in Lemma \ref{lemme_controle_energie} and \eqref{hyp} in Lemma \ref{lemme2} are fulfilled.
\end{rk}
As a consequence of Lemmas \ref{lemme_controle_energie} and \ref{lemme2}, we obtain that %Now, we shall prove that these estimations (on $\|u^\epsalpha (t) \|_{L^\infty}$ and, consequently, on $E^\delta_\epsalpha (t)$) remain true on the whole existence interval $[0,\tmax)$. 
%We have already proved that $T > 0$. Assume $T < \tmax_\eps$. Then
% \begin{equation}
% \label{contradiction}
%  \| u^\eps  (T)\|_\infty = \dfrac{1}{8 C_1 \sqrt{\eps}}.
% \end{equation}
% On the other hand, \eqref{ongardelecontrole} in Lemma \ref{lemme2} yields 
% \be
%  \|u^\eps  (T)\|_\infty \leq \dfrac{1}{16 C_1 \sqrt{\eps}} < \dfrac{1}{8 C_1 \sqrt{\eps}},
% \ee
% which contradicts \eqref{contradiction} so $T \geq \tmax_\eps$.
% We deduce then from Lemma \ref{lemme_controle_energie} and Lemma \ref{lemme2} that 
$E^\delta_\eps$ satisfies inequality \eqref{globalisation} on the whole existence interval $[0,\tmax)$. Therefore $(HNS^\epsalpha)$ has a global solution.

% Then, starting at $t=T$, we repeat the reasoning:
% $$ 
% \frac{\textrm{d}}{\textrm{d}t}E_\epsalpha^0 (T) \leq -\frac{1}{4} \|\u(T)\|_{\dot H^1}^2 + \left(\frac{C^2}{2} \alpha  \|\u(T) \|_{L^2}^2 -1 \right) \times \frac{1}{\alpha} \| \div \u(T) \|_{L^2}^2 \leq -\frac{1}{4} \|\u(T)\|_{\dot H^1}^2 
% $$ 
% since we have $$\frac{C^2}{2} \alpha  \|\u(T) \|_{L^2}^2 \leq 2 C^2 \alpha \|\u_0\|_{L^2}^2 \leq \frac{1}{2}.$$
% So $E_\epsalpha^0$ decreases in a neighborhood of $T$. We thereby obtain that $T=\tmax$ and we conclude that $E_\epsalpha^0$ decreases on $[0, \tmax).$\\

%Finally, assuming the conditions \eqref{size} on the initial data size, we obtain that $T \geq \tmax$ and therefore conclude that the local solution $\u$ to $(HNS^\epsalpha)$ is global.

% \subsection{Existence locale pour \eqref{equ}}
% On applique la m\'ethode du point fixe de Picard dans l'espace fonctionnel m\'etrique complet
% \begin{eqnarray}
%   X_T (a) & = & \left\lbrace (f,\partial_t f) \in L_T^\infty  (\dot H^{\frac{n}{2}+\delta}\cap \dot H^{\frac{n}{2}+\delta -1} )(\R^n) \times L_T^\infty   \dot H^{\frac{n}{2}+\delta -1}(\R^n) ~ : ~  \right. \nonumber \\
%  & & ~ \left. \|f\|_{X_T}~ := ~\|f\|_{L^\infty_T  \dot H^{\frac{n}{2}+\delta}} + \|f\|_{L^\infty_T \dot H^{\frac{n}{2}+\delta-1} } + \| \partial_t f \|_{L^\infty_T \dot H^{\frac{n}{2}+\delta-1}} \leq a \right\rbrace , \nonumber
% \end{eqnarray}
% avec $a>0$ et $0< T \leq 1$ \`a choisir \emph{ad hoc}.
\subsection{Convergence} \label{conv2d}
Let $\u$ and $u^\eps$ be the solutions to
\begin{eqnarray*}
 (HNS^\epsalpha)~~~  \eps \dtt \u + \dt \u - \Delta \u & = & -(\u.\nabla) \u + \frac{1}{\alpha} \nabla \div \u \\
 (HNS^\eps) ~~~~~~~~~~~ \eps \dtt u^\eps + \dt u^\eps - \Delta u^\eps & = & -(u^\eps.\nabla) u^\eps + \nabla p^\eps ~ , ~~~~~\div u^\eps =0  
\end{eqnarray*}
with the same initial data $\u_0 = u^\eps_0 \in H^{1+\delta}(\R^2)^2$ and $\u_1 = u^\eps_1 \in H^{\delta}(\R^2)^2$.\\
Let us now define the following modulated energy
$$E_{\eps ,\alpha , u^\eps} \! = \! \frac{1}{2} \|\u-u^\eps + \eps \dt (\u-u^\eps)\|_{L^2}^2 + \frac{\eps ^2}{2} \|\dt (\u-u^\eps) \|_{L^2}^2 + \eps \|\u-u^\eps \|_{\dot H^1}^2 + \frac{\eps}{\alpha} \|\div \u \|_{L^2}^2 
$$
which is inspired by the Dafermos modulated energy (see \cite{brenier, article}). Notice that
\[ E_{\eps ,\alpha , u^\eps}\left( \u \right)= E_{\eps ,\alpha}\left(\u - u^\eps \right) .\]
Through a Gronwall estimate on this energy, we shall prove that $\u$ converges to $u^\eps$ in the $L^\infty_T L^2 (\R^2)$ norm, as $\alpha$ goes to $0$. To this end, let us estimate the time derivative of the modulated energy $E_{\eps ,\alpha , u^\eps}$ using the equations $(HNS^\eps)$ and $(HNS^\epsalpha)$. We have
\begin{small}
\begin{eqnarray*}
 \frac{\d}{\d t}E_{\eps ,\alpha , u^\eps} & = & \int_{\R^2}  \bigg( \left( \partial_t (\u -u^\eps ) + \eps \partial_{tt} (\u - u^\eps) \right). \left( \u - u^\eps + \eps \partial_t (\u - u^\eps) \right)  \\
 & &  + \eps^2 \partial_{tt} (\u - u^\eps ). \partial_t (\u - u^\eps ) + 2\eps \partial_t \nabla(\u - u^\eps ).\nabla (\u - u^\eps )  \\
 & &   - \frac{2\eps}{\alpha} \partial_t \u . \nabla (\div \u ) \bigg) \, \d x \\
 &=  & \int_{\R^2} \bigg( \left( \partial_t (\u -u^\eps ) + \eps \partial_{tt} (\u - u^\eps) \right). \left( \u - u^\eps + 2\eps \partial_t (\u - u^\eps) \right)    \\
 & &  -  \eps |\partial_t (\u - u^\eps) |^2- 2\eps \partial_t (\u - u^\eps). \Delta (\u - u^\eps) - \frac{2\eps}{\alpha} \partial_t \u . \nabla (\div \u ) \bigg) \, \d x \\
%  = & \int_{\R^2} \left( \Delta (\u - u^\eps) + (u^\eps . \nabla) u^\eps - (\u . \nabla) \u + \frac{1}{\alpha} \nabla (\div \u) - \nabla p^\eps \right). \left( \u - u^\eps + 2\eps \partial_t (\u - u^\eps) \right)- \\
%   & - \eps |\partial_t (\u - u^\eps) |^2 - 2\eps \partial_t (\u - u^\eps). \Delta (\u - u^\eps) - \frac{2\eps}{\alpha} \partial_t \u . \nabla (\div \u ) \\
%  & \\
&=  & \int_{\R^2} \bigg( -|\nabla (\u - u^\eps )|^2 + \big( (u^\eps .\nabla) u^\eps - (\u . \nabla)\u \big).(\u - u^\eps )  - \eps |\partial_t (\u - u^\eps) |^2  \\
  & & + 2\eps \partial_t (\u - u^\eps ) . \big( (u^\eps .\nabla) u^\eps - (\u . \nabla)\u \big) + \frac{1}{\alpha} \nabla (\div \u ). (\u -u^\eps) + \\
  & &  + \frac{2\eps}{\alpha} \nabla (\div \u). \left( \partial_t (\u -u^\eps) - \partial_t \u \right) - \nabla p^\eps . (\u - u^\eps) - 2\eps \nabla p^\eps .\partial_t(\u - u^\eps) \bigg) \, \d x .
%  = &  \int_{\R^2} -|\nabla (\u - u^\eps )|^2 + \eps | (u^\eps .\nabla) u^\eps - (\u . \nabla)\u |^2 - \eps | \partial_t (\u - u^\eps) -(u^\eps .\nabla) u^\eps + (\u . \nabla)\u |^2 + \\
%   & + \left( (u^\eps .\nabla) u^\eps - (\u . \nabla)\u \right). (\u - u^\eps) - \frac{1}{\alpha} |\div \u |^2 -\frac{2\eps}{\alpha} \nabla (\div \u ). \partial_t u^\eps - \nabla p^\eps .(\u+ 2\eps \partial_t  \u )\\
%  & \\
\end{eqnarray*}
\end{small}
Now, using that $\div u^\eps =0$ and the inequality
\begin{small}
\[
 -\eps |\dt (\u - u^\eps)|^2 + 2\eps \dt (\u -u^\eps). \big( (u^\eps .\nabla) u^\eps - (\u . \nabla)\u \big) \leq \eps | (u^\eps .\nabla) u^\eps - (\u . \nabla)\u |^2
\]
\end{small}
and integrating in time, we finally obtain the estimate
\begin{small}
\begin{eqnarray*}
E_{\eps ,\alpha , u^\eps} (T) & \leq  & E_{\eps ,\alpha , u^\eps} (0) - \| \u - u^\eps \|_{L^2_T \dot H^1}^2 + \int_0^T \!\!\!\! \int_{\R^2} (\u - u^\eps). \left( (u^\eps .\nabla) u^\eps - (\u .\nabla) \u \right) \, \d x \, \d t +  \\
  & & + \eps \| (u^\eps .\nabla) u^\eps - (\u . \nabla)\u \|_{L^2_T L^2}^2     - \int_0^T \!\!\!\! \int_{\R^2} \nabla p^\eps .(\u + 2\eps  \partial_t \u ) \, \d x \, \d t
\end{eqnarray*}
\end{small}
In the following subsections, we shall estimate the terms in the right hand side. \\
\\
\textbf{Notation:} We shall write $f=\mathcal{O}(1)$ if there exists a constant $C = C(\eps , u^\eps)$ such that $f\leq C$. Similarly, $f=\mathcal{O} (\sqrt{\alpha})$ means that $f \leq C (\eps , u^\eps) \sqrt{\alpha}$.

\subsubsection{Estimate on $\int_0^T \int_{\R^2} (\u - u^\eps). \left( (u^\eps .\nabla) u^\eps - (\u .\nabla) \u \right) \, \d x \d t$} 
From the estimates on the energies $E^0_{\eps , \alpha}$ and $E^0_{\eps}$, we know that 
$$
 \u, u^\eps  \in  L^\infty_T L^2(\R^2) \cap L^2_T \dot H^1 (\R^2)
$$
and the boundedness of $E_\epsalpha^0$ yields
$$ \| \div \u \|_{L^2_T L^2} = \mathcal{O}(\sqrt{\alpha}).
$$
%Moreover, the Sobolev embedding $\dot H^{\frac{1}{2}}(\R^2) \subset L^4 (\R^2)$ followed by Ladyzhenskaya's inequality yields, for $f \in L^4 (\R^2)$,
Besides, let us recall Ladyzhenskaya's inequality:
\begin{equation} 
 \label{SGN} \|f\|_{L^4}^2 \leq c \|f\|_{L^2} \|f\|_{\dot H^1}
\end{equation}
which holds for all $f \in H^1(\R^2)$. 
Hence we can estimate the integral 
\[ I= \int_{\R^2} (\u - u^\eps). \left( (u^\eps .\nabla) u^\eps - (\u .\nabla) \u \right) \, \d x \]
as follows:
\begin{eqnarray*}
 I  & = & \int_{\R^2} \left[ (\u-u^\eps) . \left( (u^\eps-\u).\nabla \right) u^\eps  +  (\u-u^\eps) . (\u.\nabla) (u^\eps-\u) \right]  \, \d x \\
 & = &  \int_{\R^2} (\u-u^\eps) . \left( (u^\eps-\u).\nabla \right) u^\eps \, \d x + \frac{1}{2} \int_{\R^2} \div \u ~ |\u-u^\eps|^2  \, \d x \\
 & \leq & \|\u-u^\eps\|_{L^4}^2 \|\nabla u^\eps\|_{L^2} + \frac{1}{2} \|\div \u \|_{L^2} \| \u-u^\eps \|_{L^4}^2\\
 & \stackrel{\eqref{SGN}}{\leq} & C \| \u - u^\eps\|_{L^2} \| \u - u^\eps\|_{\dot H^1} \left( \|u^\eps \|_{\dot H^1} + \|\div \u \|_{L^2} \right) \\
 & \leq & \frac{C^2}{\eta} \left( \|u^\eps\|_{\dot H^1} + \|\div \u \|_{L^2} \right)^2 \|\u-u^\eps\|_{L^2}^2 + \eta  \|\u-u^\eps\|_{\dot H^1}^2\\
 & \leq & \frac{2 C^2}{\eta} \left( \|u^\eps\|_{\dot H^1}^2 + \|\div \u \|_{L^2}^2 \right) \| \u -u^\eps\|_{L^2}^2 + \eta  \| \u-u^\eps\|_{\dot H^1}^2 ,
\end{eqnarray*}
where $\eta>0$ is a small number to be chosen in the conclusion.

\subsubsection{Estimate on $  \eps \| (u^\eps.\nabla) u^\eps -(\u.\nabla)\u \|_{L^2_T L^2}^2 - \frac{1}{2} \|\u-u^\eps \|_{L^2_T \dot H^1}^2   $}
Let us set $$ A^{\eps , \alpha} (t) =\eps \| (u^\eps.\nabla) u^\eps -(\u.\nabla)\u \|_{L^2}^2 - \frac{1}{2} \| \u-u^\eps \|_{\dot H^1}^2 $$
and write $$ (\u.\nabla) \u - (u^\eps.\nabla) u^\eps = \big(\u.\nabla \big) (\u-u^\eps) + \big( (\u-u^\eps).\nabla \big) u^\eps $$ then, using Lemma \ref{div_lp}, estimate
\begin{eqnarray*}
 \| (u^\eps.\nabla) u^\eps \! -\! (\u.\nabla)\u \|_{L^2} & \leq &  \| \u \|_{L^\infty} \| \u -u^\eps\|_{\dot H^1} + \| \big( (\u-u^\eps).\nabla \big) u^\eps \|_{L^2 }.
\end{eqnarray*}
The second term on the right hand side estimates using the Sobolev embeddings 
$$\dot H^s \subset L^{\frac{2}{1-s}}(\R^2)$$ 
as follows. Notice that $\nabla u^\eps \in L^\infty_T \dot H^\delta \subset L^\infty_T L^{\frac{2}{1-\delta}}$ and 
$$\u - u^\eps \in L^2_T L^2 \cap L^2_T \dot H^1 \subset L^2_T \dot H^{1-\delta} \subset L^2_T L^{\frac{2}{\delta}}.$$ 
So we have
\begin{eqnarray*}
 \| \big( (\u-u^\eps).\nabla \big) u^\eps \|_{L^2} & \leq & \| \u - u^\eps \|_{L^{\frac{2}{\delta}}} \| \nabla u^\eps \|_{L^{\frac{2}{1-\delta}}} \\
 & \leq & C \| \u - u^\eps \|_{L^2}^{\delta} \| \u - u^\eps \|_{\dot H^1}^{1-\delta} \| \nabla u^\eps \|_{\dot H^\delta}  .
\end{eqnarray*}
Now, integrating in time the squared norm and performing a Young inequality, we obtain
\begin{eqnarray*}
  \| \big (\u \!\! -u^\eps).\nabla \big) u^\eps \|_{L^2_T L^2}^2 & \!\!\!\leq &\!\!\! \frac{C_\delta}{\eta}\!\!  \int_0^T \!\!\! \|u^\eps \|_{\dot H^{1+\delta}}^2 E_{\epsalpha , u^\eps}(t) \, \d t + \tilde C_\delta \eta \|u^\eps \|_{L^\infty_T \dot H^{1+\delta}}^2 \|\u-u^\eps \|_{L^2_T \dot H^1}^2.
\end{eqnarray*}
Finally, recalling that $\eps \left( \| u^\eps \|_{ L^\infty_T \dot H^{1+\delta} }^2 + \| u^\eps \|_{L^2_T \dot H^{1+\delta}}^2 \right) =\mathcal{O}(1) $, we obtain the inequality
\begin{eqnarray*}
 \int_0^T \!\!\! A^{\eps , \alpha} \d t & \leq & \left( C \eps \| \u\|_{L^\infty_T L^\infty}^2 + \tilde C_\delta \eps \eta \|u^\eps\|_{L^\infty_T \dot H^{1+\delta}}^2- \frac{1}{2} \right) \| \u-u^\eps\|_{L^2_T \dot H^1}^2 + \\
 &      & + \frac{C_\delta}{\eta} \int_0^T \!\!\! \eps \| u^\eps\|_{\dot H^{1+\delta}}^2 E_{\epsalpha , u^\eps}(t) \d t\\
 & \leq & \frac{C_\delta}{\eta} \int_0^T \!\!\! \eps \| u^\eps\|_{\dot H^{1+\delta}}^2 E_{\epsalpha , u^\eps}(t) \d t
\end{eqnarray*}
if $\eps$ and $\eta$ are small enough.

\subsubsection{Estimate on $\int_0^T \int_{\R^2} \nabla p^\eps .\u \, \d x \, \d t$}
First, notice that, applying the $\div $ operator to equation $(HNS^\eps)$, we obtain the identity
\begin{equation} \label{pv} \Delta p^\eps = \div (u^\eps.\nabla) u^\eps = \div \nabla : u^\eps \otimes u^\eps \end{equation}
from which we deduce that $p^\eps$ has the same regularity as $u^\eps \otimes u^\eps$ and 
$$\|p^\eps \|_{L^2_T L^2} \leq C \| u^\eps \otimes u^\eps \|_{L^2_T L^2} .$$ 
Since $u^\eps \in L^2_T \dot H^1 \cap L^\infty_T L^2 \subset L^4_T L^4 (\R^2)^2$, we immediately conclude that 
\begin{eqnarray*}
 \int_0^T \int_{\R^2} \nabla p^\eps . \u  \, \d x \, \d t & = & -  \int_0^T \int_{\R^2} p^\eps . \div \u  \, \d x \, \d t 
 ~ \leq ~ \|p^\eps \|_{L^2_T L^2} \| \div \u \|_{L^2_T L^2} \\
 & \leq & C \sqrt{\alpha} \| p^\eps \|_{L^2_T L^2} 
 ~ \leq ~ C \sqrt{\alpha} \| u^\eps \otimes u^\eps \|_{L^2_T L^2} \\
 & \leq & C \sqrt{\alpha} \|u^\eps\|_{L^4_T L^4}^2 ~=~ \mathcal{O}(\sqrt{\alpha}) .
\end{eqnarray*}
\subsubsection{Estimate on $2 \eps \int_0^T \int_{\R^2} \nabla p^\eps . \partial_t \u \, \d x \, \d t$}
From identity \eqref{pv}, we know that $$ \partial_t p^\eps = \sum_{i,j=1}^2 \partial_t \frac{\partial_i \partial_j}{\Delta} (u^\eps_i u^\eps_j) .$$
Besides, due to the control of the energies $E_\eps^0$ and $E^\delta_\eps$, we have
\begin{eqnarray*}
 \partial_t u^\eps & \in & L^2_T L^2 \\%\cap L^\infty_T H^s \\
 u^\eps & \in & L^\infty_T H^{1+\delta} \subset L^\infty_T L^\infty . %L^\infty_T L^2 \cap L^2_T H^1 \subset L^4_T L^4 \\
\end{eqnarray*}
After two integrations by parts (one in each variable), we obtain an integral which easily estimates:
\begin{eqnarray*}
 \int_0^T \int_{\R^2} \nabla p^\eps . \partial_t \u \, \d x \, \d t & = &  \int_0^T \int_{\R^2} \partial_t p^\eps . \div \u \, \d x \, \d t \\
 & \leq & \| \partial_t p^\eps \|_{L^2_T L^2} \| \div \u \|_{L^2_T L^2} \\
 & \leq & C \sqrt{\alpha} \| \partial_t u^\eps \|_{L^2_T L^2} \| u^\eps \|_{L^\infty_T L^\infty}  ~=~ \mathcal{O}(\sqrt{\alpha}) .
\end{eqnarray*}

\subsubsection{Conclusion}
First, notice that since we take the same initial data for $(HNS^\epsalpha)$ and $(HNS^\eps)$, we have in particular $\div \u_0 =0$ so
$$E_{\epsalpha , u^\eps}(0)=0.$$ 
Now, gathering the estimates in the previous subsections, we obtain that
\begin{eqnarray*}
E_{\epsalpha , u^\eps}(T) &\leq & \mathcal{O}(\sqrt{\alpha}) + C_{\delta , \eta}\int_0^T \left( \|u^\eps \|_{\dot H^1}^2 + \eps \|u^\eps \|_{\dot H^{1+\delta}}^2 + \mathcal{O}(\alpha) \right) E_{\epsalpha , u^\eps}(t) \, \d t \\
& & + \left( \eta + (1-\delta) \eps \eta \|u^\eps\|_{L^\infty_T \dot H^{1+\delta}}^2 -\frac{1}{2} \right) \| \u - u^\eps\|_{L^2_T \dot H^1}^2   .
\end{eqnarray*}
So, choosing $\eta$ small enough, Gronwall's lemma yields, for all positive $T$,
$$ E_{\epsalpha , u^\eps}(T) = \mathcal{O}(\sqrt{\alpha}) .$$
Now, recall that Theorem \ref{th} tells that $u^\eps$ converges towards the solution $v$ to $(NS)$. Theorem \ref{th_alpha} is now proved.
\section{The 3D case} \label{3d}
\subsection{Preliminary estimates} \label{rappels3d}           
As in subsection \ref{rappels2d}, we start by recalling the regularity results we have on $u^\eps$. All the proofs are in \cite{article}.
\begin{lemme} Let $u^\eps \in L^\infty_T (\dot H^{\frac{3}{2}+\delta} \cap \dot H^{\frac{1}{2}+\delta})(\R^3)^3$ be the global solution to $(HNS^\eps)$ with initial data $(u_0^\eps , u^\eps_1 ) \in H^{\frac{3}{2}+\delta} \times H^{\frac{1}{2}+\delta} (\R^3)^3$. Then we have
 $$u^\eps \in L^\infty_T \dot H^{\frac{3}{2}+\delta} \cap L^2_T \dot H^{\frac{3}{2}+\delta} \cap L^2_T \dot H^{\frac{3}{2}} \cap L^\infty_T \dot H^{\frac{1}{2}} \cap L^\infty_T \dot H^{\frac{3}{2}}$$
 and
 $$\dt u^\eps \in L^2_T \dot H^{\frac{1}{2}} \cap L^\infty_T \dot H^{\frac{1}{2}+\delta}.$$
\end{lemme}
\begin{rk}
In this paper, we are not interested in the dependence of the norms on $\eps$.
\end{rk}
\proof Here is a sketch of the proof. First, let us define the 3D energy
$$ E_\eps^{\frac{1}{2}+\delta} (t) = \frac{1}{2} \| u^\eps + \eps \dt u^\eps \|_{\dot H^{\frac{1}{2}+\delta}}^2 + \frac{\eps^2}{2} \| \dt u^\eps \|_{\dot H^{\frac{1}{2}+\delta}}^2 + \eps \| \nabla u^\eps \|_{\dot H^{\frac{1}{2}+\delta}}^2
$$ for non-negative $\delta$. In order to obtain that the solutions $u^\eps$ to $(HNS^\eps)$ are global, we proved in \cite{article} that
\be \label{globalisation} \exists C>0 ~:~ \forall t \geq 0 ~,~~ E_\eps^{\frac{1}{2}+\delta} (t) \leq C \eps^{-\delta} . \ee
Directly from the expression of $E^\delta_\eps$ and from \eqref{globalisation}, we deduce that 
$$\partial_t u^\varepsilon \in L^\infty_T  \dot H^{\frac{1}{2}+\delta} ~,~~    u^\eps \in L^\infty_T \dot H^{\frac{3}{2}+\delta} \cap L^2_T \dot H^{\frac{3}{2}+\delta}. $$
Besides, let us consider the time derivative of $E_\eps^{\frac{1}{2}}$. According to \cite{article}, we have
$$
E_\eps^{\frac{1}{2}}(T) + \eps \| \dt u^\eps + (u^\eps . \nabla)u^\eps \|_{L^2_T \dot H^{\frac{1}{2}}}^2 + \|u^\eps \|_{L^2_T \dot H^{\frac{3}{2}}}^2 - \eps \| (u^\eps . \nabla)u^\eps \|_{L^2_T \dot H^{\frac{1}{2}}}^2 = E_\eps^{\frac{1}{2}}(0) + I(\eps , T) , 
$$
where 
\begin{eqnarray*}
 I(\eps , T) & := &\int_0^T \!\!\! \int_{\R^3} \Lambda^{\frac{1}{2}} \left( (u^\eps .\nabla) u^\eps \right). \Lambda^{\frac{1}{2}} u^\eps \, \d x \, \d t \\
 & \leq & \|u^\eps \|_{L^\infty_T \dot H^{\frac{1}{2}}} \, \|u^\eps \|_{L^2_T \dot H^{\frac{3}{2}}}^2 .
\end{eqnarray*}
The $\Lambda$ above is the Fourier multiplier defined by $\widehat{\Lambda f}(\xi) =|\xi| \hat f(\xi).$\\
So, since $E_\eps^{\frac{1}{2}}(0) \leq C$ and $\|u^\eps \|_{L^\infty_T \dot H^{\frac{1}{2}}} \leq 2 \|u^\eps_0\|_{\dot H^{\frac{1}{2}}} \leq \frac{1}{8} $ by the decay of $E_\eps^{\frac{1}{2}}$, we obtain the inequality
$$ E_\eps^{\frac{1}{2}}(T) + \eps \| \dt u^\eps + (u^\eps . \nabla)u^\eps \|_{L^2_T \dot H^{\frac{1}{2}}}^2 + \frac{1}{2} \|u^\eps \|_{L^2_T \dot H^{\frac{3}{2}}}^2 \leq C $$
if $\eps$ is small enough. Now, we complete the reasoning as in the 2D case and deduce that
$$ \dt u^\eps \in L^2_T \dot H^{\frac{1}{2}} ~,~~ u^\eps \in L^2_T \dot H^{\frac{3}{2}} . $$
Finally, since $E_\eps^{\frac{1}{2}}$ is bounded, we immediately have
$$u^\eps \in L^\infty_T \dot H^{\frac{1}{2}} \cap L^\infty_T \dot H^{\frac{3}{2}}$$
\subsection{Globalization} \label{glob3d}
As in the 2D case, we shall prove that the energy
$$E_{\eps , \alpha}^{\frac{1}{2}+\delta}(t) =  \frac{1}{2} \|\u + \eps \partial_t \u  \|_{\dot H^{\frac{1}{2}+\delta}}^2 + \frac{\eps^2}{2} \|\partial_t \u \|_{\dot H^{\frac{1}{2}+\delta} }^2 + \eps \|\u \|_{\dot H^{\frac{3}{2}+\delta}}^2 + \frac{\eps}{\alpha} \|\div \u \|_{\dot H^{\frac{1}{2}+\delta} }^2 
$$
is bounded in order to show that the local solution obtained by the contraction argument is global. First, let us define $0 \leq T \leq \tmax $ by
\be
 T = \sup \left\lbrace 0 \leq \tau \leq \tmax ~ :~ \forall ~ t \in [0, \tau) , ~\|\u  (t)\|_{L^\infty} < \dfrac{1}{2 K_1 \sqrt{\eps}} \right\rbrace ,
\ee
where $K_1=\max (K, \tilde K)$ such that
\begin{eqnarray} 
\|(\u.\nabla) \u\|_{\dot H^{\frac{1}{2}}} & \leq & K \|\u\|_{L^\infty} \|\u\|_{\dot H^{\frac{3}{2}}} , \label{nonlinearite3D} \\
\|(\u.\nabla) \u\|_{\dot H^{\frac{1}{2}+\delta}} & \leq & \tilde K \|\u\|_{L^\infty} \|\u\|_{\dot H^{\frac{3}{2}+\delta}} . \label{nonlinearite3Ddelta} 
\end{eqnarray}
We shall start by proving that the energy $E^{\frac{1}{2}}_\epsalpha$ decreases on $[0,T)$. To this end, we estimate its time derivative:
\begin{eqnarray*}
 \frac{\d E_{\eps, \alpha}^{\frac{1}{2}}}{\d t}& \leq &  \eps \|(\u.\nabla) \u\|_{\dot H^{\frac{1}{2}}}^2  - \|\u\|_{\dot H^{\frac{3}{2}}}^2 - \int_{\R^3}  \Lambda^{\frac{1}{2}} (\u.\nabla)\u . \Lambda^{\frac{1}{2}} \u \, \d x    \\
 &   \leq &  \left( K_1^2 \eps \|\u\|_{L^\infty}^2 -1 \right) \|\u\|_{\dot H^{\frac{3}{2}}}^2  -  \int_{\R^3}  \Lambda^{\frac{1}{2}} (\u.\nabla)\u . \Lambda^{\frac{1}{2}} \u \, \d x   .
\end{eqnarray*} 
 So, for $\eps$ small enough, we have
 \[
   \frac{\d E_{\eps, \alpha}^{\frac{1}{2}}}{\d t}  \leq   - \frac{3}{4} \|\u\|_{\dot H^{\frac{3}{2}}}^2 + \int_{\R^3}  (\u.\nabla)\u . \Lambda \u \, \d x  \]
on $[0,T)$. Then, recalling that $\div \u \neq 0$, we have
\begin{eqnarray*}
 \int_{\R^3} (\u.\nabla)\u . \Lambda \u \, \d x & = & \sum_{i,j=1}^3 \int_{\R^3} \u_i \, \partial_i \u_j \, \partial_j \u_i \, \d x \\
 & = & -\int_{\R^3} \left( \sum_{i,j=1}^3  \partial_j \u_i \, \partial_i \u_j \, \u_i  + \sum_{i=1}^3 (\u_i) ^2 \, \partial_i \div \u \right) \, \d x  \\
 & = &   \sum_{i=1}^3 \int_{\R^3} \u_i \, \partial_i \u_i \ \div \u \d x \\
 & \leq &   \sum_{i=1}^3  \| \u_i \|_{L^3} \, \| \partial_i \u_i \|_{L^3} \, \| \div \u \|_{L^3} \\
 & \leq & 3 K_2^3 \| \u \|_{\dot H^{\frac{1}{2}}} \, \| \u \|_{\dot H^{\frac{3}{2}}} \, \| \div \u \|_{\dot H^{\frac{1}{2}}} \\
 & \leq &  9 K_2^3 \| \u \|_{\dot H^{\frac{1}{2}}} \, \| \u \|_{\dot H^{\frac{3}{2}}}^2 ,
\end{eqnarray*}
where $K_2$ is the constant such that
$$\|f\|_{L^3(\R^3)} \leq  K_2 \|f\|_{\dot H^{\frac{1}{2}}(\R^3)} .$$
Now, assume that
\begin{equation} \label{normepetite} \|\u_0\|_{\dot H^{\frac{1}{2}}} \leq \frac{1}{36 K_2^3}. \end{equation} 
Then we obtain
$$ \frac{\d}{\d t}E_{\eps, \alpha}^{\frac{1}{2}}(0) \leq \left( - \frac{3}{4} + 9 K_2^3 \|\u_0\|_{\dot H^{\frac{1}{2}}} \right) \| \u_0 \|_{\dot H^{\frac{3}{2}}}^2 < -\frac{1}{4} \|\u_0\|_{\dot H^{\frac{3}{2}}}^2 $$
and $E_{\eps, \alpha}^{\frac{1}{2}}$ therefore decreases on an interval $[0,\tau]$. Set
$$ \tau = \sup \left\lbrace t < T ~:~ E_{\epsalpha}^{\frac{1}{2}}  \textrm{ decreases on } [0,t] \right\rbrace. $$
Since $E_{\eps, \alpha}^{\frac{1}{2}}$ decreases on $[0,\tau]$, we have 
$$ \| \u(\tau)\|_{\dot H^{\frac{1}{2}}}^2 \leq 2 E_{\eps, \alpha}^{\frac{1}{2}}(\tau) \leq E_{\eps, \alpha}^{\frac{1}{2}}(0) \leq 4 \|\u_0\|_{\dot H^{\frac{1}{2}}}^2 \leq \left( \frac{1}{18K_2^3} \right)^2 $$
so that
$$  \frac{\d}{\d t}E_{\eps, \alpha}^{\frac{1}{2}}(\tau)  \leq -\frac{1}{4} \| \u(\tau)\|_{\dot H^{\frac{3}{2}}}^2.$$
Thus $E_{\eps, \alpha}^{\frac{1}{2}}$ decreases on $[0,\tau ']$, where $\tau' > \tau$. We have proved by contradiction that the energy decays on the whole interval $[0,T)$.\\
\\
Now, we can prove the following lemma:
\begin{lemme}
 Assume that $\|\u_0\|_{\dot H^{\frac{1}{2}}} <\dfrac{1}{36 K_2^3} $ and
\be
 (H') ~~~~~
 i) ~ \eps^{\frac{1+\delta}{2}} \|\u_0 \|_{\dot{H}^{\frac{3}{2}+\delta}}  =o(1)  ~,~~ ~~~~
 ii)~ \eps^{\frac{\delta}{2}} \|\u_0 \|_{\dot{H}^{\frac{1}{2}+\delta}}  =o(1) .
\ee
when $\eps$ goes to zero. \\
Now, recall that $0 \leq T \leq \tmax $ is defined by
\[
 T = \sup \left\lbrace 0 \leq \tau \leq \tmax ~ :~ \forall ~ t \in [0, \tau) , ~\|\u  (t)\|_{L^\infty} < \dfrac{1}{2K_1 \sqrt{\eps}} \right\rbrace .
\]
Then $T>0$ and there exists a large number $N$, depending only on $\delta$, and a constant $C>1$ such that 
\be 
 E^{\frac{1}{2}+\delta}_\eps (t) \leq C^N \, E^{\frac{1}{2}+\delta}_\eps (0) 
\ee
for all $t\in [0,T)$ and $\eps$ small enough.
\end{lemme}
\proof First, let us compute the time derivative of this energy:
\newcommand{\ddt}[1]{\frac{d}{d t}#1}
\begin{eqnarray*}
 \ddt{E_{\eps , \alpha}^{\frac{1}{2}+\delta}} (t) %& = & \int_{\R^3} \Lambda^{\frac{1}{2}+\delta} (\dt u + \eps \dtt u).\Lambda^{\frac{1}{2}+\delta}(u+\eps \dt u ) + \eps^2 \Lambda^{\frac{1}{2}+\delta} \dt u . \Lambda^{\frac{1}{2}+\delta} \dtt u +\\
% & & \hfill +2\eps \Lambda^{\frac{3}{2}+\delta} \dt u .\Lambda^{\frac{3}{2}+\delta} u + \frac{2\eps}{\alpha} \Lambda^{\frac{1}{2}+\delta} \div(\dt u) . \Lambda^{\frac{1}{2}+\delta} \div u \\
 & = & \eps \|(\u.\nabla)\u \|_{\dot H^{\frac{1}{2}+\delta}}^2- \|\u\|_{\dot H^{\frac{3}{2}+\delta}}^2 - \int_{\R^3} \Lambda^{\frac{1}{2}+\delta} (\u.\nabla)\u .\Lambda^{\frac{1}{2}+\delta} \u \, \d x  -  \\
 & & - \eps \| \dt \u + (\u.\nabla)\u \|_{\dot H^{\frac{1}{2}+\delta}}^2 - \frac{1}{\alpha} \| \div \u\|_{\dot H^{\frac{1}{2}+\delta}}^2 \\
 & \leq & \eps \|(\u.\nabla)\u \|_{\dot H^{\frac{1}{2}+\delta}}^2 - \|\u\|_{\dot H^{\frac{3}{2}+\delta}}^2 - \int_{\R^3} \Lambda^{\frac{1}{2}+\delta} (\u.\nabla)\u .\Lambda^{\frac{1}{2}+\delta} \u \, \d x \\
 & \leq &  \left( K_1^2 \eps \|\u\|_{L^\infty}^2 -1 \right) \|\u\|_{\dot H^{\frac{3}{2}+\delta}}^2 - \int_{\R^3} \Lambda^{\frac{1}{2}+\delta} (\u.\nabla)\u .\Lambda^{\frac{1}{2}+\delta} \u \, \d x .
 \end{eqnarray*}
 On the time interval $[0,T)$, we have $K_1^2 \eps \|\u\|_{L^\infty}^2 -1 \leq -\frac{1}{2}$ and we can estimate the integral on the right hand side by
 \begin{eqnarray*}
  \int_{\R^3} \Lambda^{\frac{1}{2}+\delta} (\u.\nabla)\u .\Lambda^{\frac{1}{2}+\delta} \u \, \d x & \leq & \|(\u.\nabla)\u \|_{\dot H^{\frac{1}{2}+\delta}} \|\u \|_{\dot H^{\frac{1}{2}+\delta}} \\
  & \leq & K_1 \|\u\|_{L^\infty} \|\u\|_{\dot H^{\frac{3}{2}+\delta}} \|\u \|_{\dot H^{\frac{1}{2}+\delta}}
 \end{eqnarray*}
due to Lemma \ref{div_lp}. Using that followed by \eqref{interp_besov}, we obtain 
  \begin{eqnarray*}
\ddt{E_{\eps , \alpha}^{\frac{1}{2}+\delta}} (t) %& \leq & - \frac{1}{2} \|\u\|_{\dot H^{\frac{3}{2}+\delta}}^2 + \|(\u.\nabla)\u \|_{\dot H^{\frac{1}{2}+\delta}} \|\u \|_{\dot H^{\frac{1}{2}+\delta}} \\
 & \leq & - \frac{1}{2} \|\u\|_{\dot H^{\frac{3}{2}+\delta}}^2 + K_1 \|\u\|_{L^\infty} \|\u\|_{\dot H^{\frac{3}{2}+\delta}} \|\u \|_{\dot H^{\frac{1}{2}+\delta}} \\
 & \leq & - \frac{1}{2} \|\u\|_{\dot H^{\frac{3}{2}+\delta}}^2 + C \|\u \|_{\dot H^{\frac{1}{2}+\delta}}^{1+\delta} \|\u\|_{\dot H^{\frac{3}{2}+\delta}}^{2-\delta} \\
 & \leq & - \frac{1}{2} \|\u\|_{\dot H^{\frac{3}{2}+\delta}}^2 + C \|\u \|_{\dot H^{\frac{1}{2}}}^{1-\delta} \|\u \|_{\dot H^{\frac{3}{2}}}^{\delta} \|\u \|_{\dot H^{\frac{1}{2}+\delta}}^{\delta} \|\u\|_{\dot H^{\frac{3}{2}+\delta}}^{2-\delta} \\
 & \leq & C_ \delta \|\u \|_{\dot H^{\frac{1}{2}}}^{2\frac{1-\delta}{\delta}} \|\u \|_{\dot H^{\frac{3}{2}}}^2 E_{\eps , \alpha}^{\frac{1}{2}+\delta}(t),
\end{eqnarray*}
where we have used standard interpolations and a Young inequality. Now, since we do not know if the energy $E_{\eps , \alpha}^{\frac{1}{2}+\delta}$ decays, we will use that $E_{\eps, \alpha}^{\frac{1}{2}}$ satisfies the inequality
\[ \ddt{E_\epsalpha^{\frac{1}{2}}}(t) \leq - \frac{1}{4} \|\u\|_{\dot H^{\frac{3}{2}}}^2 
\]
on the interval $[0,T)$.\\
Now, as in the 2D case, let us define the functional $\mathcal{E}_{\eps , \alpha , N}^{\frac{1}{2}} := E_{\eps, \alpha}^{\frac{1}{2}+\delta} \, \left( 1+ E_{\eps, \alpha}^{\frac{1}{2}} \right)^N$. Then we have
\begin{eqnarray*}
 \frac{\d}{\d t} \mathcal{E}_{\eps , \alpha , N}^{\frac{1}{2}} (t) & = &  \frac{\d}{\d t}E_{\eps, \alpha}^{\frac{1}{2}+\delta} \, \left( 1+ E_{\eps, \alpha}^{\frac{1}{2}} \right)^N + N E_{\eps, \alpha}^{\frac{1}{2}+\delta} \frac{\d}{\d t}E_{\eps, \alpha}^{\frac{1}{2}} \, \left( 1+ E_{\eps, \alpha}^{\frac{1}{2}} \right)^{N-1} \\
 & \leq & E_{\eps, \alpha}^{\frac{1}{2}+\delta}  \left( 1+ E_{\eps, \alpha}^{\frac{1}{2}} \right)^{N-1}  \!\!\!\!\! \| \u(t) \|_{\dot H^{\frac{3}{2}}}^2 \, \left[ C_\delta \| \u(t) \|_{\dot H^{\frac{1}{2}}}^{2 \frac{1-\delta}{\delta}} \left( 1+ E_{\eps, \alpha}^{\frac{1}{2}} \right) - \frac{N}{4} \right] \\
 & \leq & E_{\eps, \alpha}^{\frac{1}{2}+\delta}  \left( 1+ E_{\eps, \alpha}^{\frac{1}{2}} \right)^{N-1} \!\!\!\!\! \| \u(t) \|_{\dot H^{\frac{3}{2}}}^2 \, \left[ C_\delta \| \u_0 \|_{\dot H^{\frac{1}{2}}}^{2 \frac{1-\delta}{\delta}} \left( 1+ 2\| \u_0 \|_{\dot H^{\frac{1}{2}}}^2 \right) - \frac{N}{4} \right] .
\end{eqnarray*}
Taking $N=N(\delta)$ large enough, we obtain that $\mathcal{E}_{\eps , \alpha , N}^{\frac{1}{2}}$ decays and we have
$$ E_{\eps, \alpha}^{\frac{1}{2}+\delta} (t) \leq E_{\eps, \alpha}^{\frac{1}{2}+\delta}(0) \, \left( 1+ E_{\eps, \alpha}^{\frac{1}{2}}(0) \right)^N. $$
Finally, notice that \[
                      E_{\eps, \alpha}^{\frac{1}{2}}(0) \leq 2 \|\u_0\|_{\dot H^{\frac{1}{2}}}^2,
                     \]
so the proof of the lemma is finished. \qed \\

As in the 2D case, the smallness assumptions on the initial data yield the boundedness of $E_{\eps, \alpha}^{\frac{1}{2}+\delta}$ for $\alpha$ and $\eps$ small enough. From this, we deduce that the solutions to $(HNS^\epsalpha)$ are global.
\subsection{Convergence} \label{conv3d}
Let $\u$ and $u^\eps$ be the solutions to $(HNS^\epsalpha)$ and $(HNS^\eps)$ respectively with the same initial data 
\[ (u_0^{\eps , \alpha} , u_1^{\eps , \alpha})= ( u^\eps_0 , u^\eps_1) \in H^{\frac{3}{2}+\delta}(\R^3)^3\times H^{\frac{1}{2}+\delta}(\R^3)^3  .\]
In order to prove that $\u$ converges to $u^\eps$ in the $L^\infty_T \dot H^{\frac{1}{2}}$ norm, as $\alpha$ goes to $0$, we shall use, as in the 2D case, a variant of the Dafermos modulated energy:
\begin{small}
\[ \mathcal{E}_{\eps , \alpha ,u^\eps}(t) :=  \frac{1}{2} \|\u-u^\eps + \eps \dt (\u-u^\eps) \|_{\dot H^{\frac{1}{2}}}^2 + \frac{\eps^2}{2} \|\dt (\u-u^\eps) \|_{\dot H^{\frac{1}{2}}}^2 + \eps \|\u-u^\eps \|_{\dot H^{\frac{3}{2}}}^2 + \frac{\eps}{\alpha} \| \div \u \|_{\dot H^{\frac{1}{2}}}^2 . \]
\end{small}
As done in section \ref{conv2d}, we shall compute the time derivative of $\mathcal{E}_{\eps , \alpha ,u^\eps}$ and use equations $(HNS^\epsalpha)$ and $(HNS^\eps)$. Doing so, we obtain the following estimate:
\begin{eqnarray}
 \mathcal{E}_{\eps , \alpha ,u^\eps}(T) & \leq & \mathcal{E}_{\eps , \alpha ,u^\eps}(0) + \eps \| (\u.\nabla)\u - (u^\eps.\nabla)u^\eps \|_{L^2_T \dot H^{\frac{1}{2}}}^2 - \| \u-u^\eps  \|_{L^2_T \dot H^{\frac{3}{2}}}^2  - \nonumber \\
 & &  - \int_0^T \int_{\R^3} \Lambda^{\frac{1}{2}}p^\eps . \Lambda^{\frac{1}{2}} \div \u \, \d x \d t + 2\eps \int_0^T \!\!\! \int_{\R^3} \Lambda^{\frac{1}{2}} \dt p^\eps . \Lambda^{\frac{1}{2}} \div \u \, \d x \d t + \nonumber \\ 
 & & + \int_0^T \!\!\!\!\! \int_{\R^3} \left[ (\u.\nabla) \u - (u^\eps.\nabla)u^\eps \right] . \Lambda (\u-u^\eps) \, \d x \, \d t- \frac{1}{\alpha} \| \div \u  \|_{L^2_T \dot H^{\frac{1}{2}}}^2 . \nonumber
\end{eqnarray}

\subsubsection{Estimate on $ \eps \| (\u.\nabla)\u - (u^\eps.\nabla)u^\eps \|_{L^2_T \dot H^{\frac{1}{2}}}^2$}
First, $A:=(\u.\nabla)\u - (u^\eps.\nabla)u^\eps$ writes
$$ (\u.\nabla)\u - (u^\eps.\nabla)u^\eps = \big((\u-u^\eps).\nabla \big) u^\eps + \big( \u.\nabla \big) (\u-u^\eps) . $$
Let us estimate the RHS using the dyadic Littlewood-Paley decomposition. First, by Lemma \ref{paraproduct} in the appendix, we have
\begin{eqnarray*}
\| (\u-u^\eps).\nabla u^\eps \|_{\dot H^{\frac{1}{2}}}^2 &\leq& \sum_{i,j=1}^3 \sum_{p \in \mathbb{Z}} 2^p \| \Delta_p (\u_i - u^\eps_i ) \, . \, S_{p+1} (\partial_i u^\eps_j) \|_{L^2}^2 \\
&+& \sum_{i,j=1}^3 \sum_{q \in \mathbb{Z}} 2^q \| \Delta_q (\partial_i u^\eps_j) \, . \, S_q (\u_i - u^\eps_i ) \|_{L^2}^2.
\end{eqnarray*}
The first term estimates easily:
\begin{eqnarray*}
 \sum_p 2^p \| \Delta_p (\u_i - u^\eps_i ) \, . \, S_{p+1} (\partial_i u^\eps_j) \|_{L^2}^2 & \leq & \sum_p 2^p  \| \Delta_p (\u_i - u^\eps_i )\|^2_{L^2} \, \| S_{p+1} (\partial_i u^\eps_j) \|^2_{L^\infty} \\
 & \leq & \|u^\eps \|_{L^\infty}^2 \sum_p 2^{3p} \, \| \Delta_p (\u_i -u^\eps_i) \|_{L^2}^2   \\
 & \leq & \|u^\eps \|_{L^\infty}^2  \|\u -u^\eps \|_{\dot H^{\frac{3}{2}}}^2 
\end{eqnarray*}
due to property \eqref{equivalence_des_normes} in the appendix. Now, notice that Sobolev embeddings imply that $u^\eps \in \dot H^{\frac{3}{2}+\delta} \subset W^{\frac{3}{2}, \alpha}(\R^3)$, where $\alpha = \frac{6}{3-2\delta}$. Then, if $\bar \alpha = \frac{3}{\delta}$, we have
\begin{eqnarray*}
 \sum_q 2^q \| \Delta_q (\partial_i u^\eps_j) \, . \, S_q (\u_i - u^\eps_i ) \|_{L^2}^2 & \leq & \sum_q 2^q \| \Delta_q (\partial_i u^\eps_j)\|_{L^\alpha}^2 \, \| S_q (\u_i - u^\eps_i ) \|_{L^{\bar \alpha}}^2 \\
 & \leq & \| \u - u^\eps  \|_{L^{\bar \alpha}}^2 \sum_q 2^{q+3q \left( 1-\frac{2}{\alpha} \right)} \| \Delta_q (\partial_i u^\eps_j)\|_{L^2}^2   \\
 & \leq & C \| \u - u^\eps  \|_{\dot H^{\frac{3}{2}-\delta}}^2 \sum_q 2^{3q+3q \left( 1-\frac{2}{\alpha} \right)} \| \Delta_q u^\eps_j\|_{L^2}^2   \\
 & = & C \| \u - u^\eps  \|_{\dot H^{\frac{3}{2}-\delta}}^2 \sum_q 2^{6q - \frac{6q}{ \alpha} } \| \Delta_q u^\eps_j\|_{L^2}^2   \\
 & = & C \| \u - u^\eps  \|_{\dot H^{\frac{3}{2}-\delta}}^2 \sum_q 2^{2q \left( \frac{3}{2}+\delta \right) } \| \Delta_q u^\eps_j\|_{L^2}^2   \\
 & \leq & C  \| \u - u^\eps  \|_{\dot H^{\frac{1}{2}}}^\delta \, \| \u - u^\eps  \|_{\dot H^{\frac{3}{2}}}^{1-\delta} \| u^\eps \|_{\dot H^{\frac{3}{2}+\delta}}  \\
 & \leq & C_\delta  \left( \| \u - u^\eps  \|_{\dot H^{\frac{1}{2}}} + \| \u - u^\eps  \|_{\dot H^{\frac{3}{2}}} \right) \| u^\eps\|_{\dot H^{\frac{3}{2}+\delta}} . \\
\end{eqnarray*}
Now, let us estimate 
\begin{eqnarray*}
\|\u.\nabla(\u-u^\eps) \|_{\dot H^{\frac{1}{2}}} &\leq &\sum_{i,j=1}^3 \sum_{p \in \mathbb{Z}} 2^p \| \Delta_p \u_i  \, . \, S_{p+1} \partial_i (\u_j - u^\eps_j) \|_{L^2}^2 \\
&+& \sum_{i,j=1}^3 \sum_{q \in \mathbb{Z}} 2^q \| \Delta_q \partial_i (\u_j - u^\eps_j) \, . \, S_q \u_i  \|_{L^2}^2.
\end{eqnarray*}
We know that $\u \in \dot H^{\frac{3}{2}+\delta} \subset \dot W^{\frac{3}{2}, \alpha}$, where $\alpha = \frac{6}{3-2\delta}$. Let $\bar \alpha = \frac{3}{\delta}$. Then
\begin{eqnarray*}
\sum_{p} 2^p  \| \Delta_p \u_i  \, . \, S_{p+1} \partial_i (\u_j - u^\eps_j) \|_{L^2}^2 & \leq & \sum_{p} 2^p  \| \Delta_p \u_i  \|_{L^\alpha}^2    \| S_{p+1} \partial_i (\u_j - u^\eps_j) \|_{L^{\bar \alpha}}^2 \\
 & \leq & \sum_{p} 2^{3p}  \| \Delta_p \u_i  \|_{L^\alpha}^2    \| S_{p+1}  (\u_j - u^\eps_j) \|_{L^{\bar \alpha}}^2 \\ 
 & \leq &  \| \u - u^\eps \|_{L^{\bar \alpha}}^2 \sum_{p} 2^{3p}  \| \Delta_p \u_i  \|_{L^\alpha}^2  \\ 
 & \leq &  \| \u - u^\eps \|_{L^{\bar \alpha}}^2 \sum_{p} 2^{3p+6p\left( \frac{1}{2}- \frac{1}{\alpha} \right)}  \| \Delta_p \u_i  \|_{L^2}^2  \\ 
 & = &  \| \u - u^\eps \|_{L^{\bar \alpha}}^2 \sum_{p} 2^{2p\left( \frac{3}{2}+\delta \right)}  \| \Delta_p \u_i  \|_{L^2}^2  \\
 & \leq &  C \| \u - u^\eps \|_{\dot H^{\frac{3}{2}-\delta}}^2  \|  \u  \|_{\dot H^{\frac{3}{2}+\delta}}^2  
\end{eqnarray*}
since $\dot H^{\frac{3}{2}-\delta} \subset L^{\bar \alpha} (\R^3)$. So, by interpolation and a Young inequality, we get
$$
 \sum_{p} 2^p  \| \Delta_p \u_i   .  S_{p+1} \partial_i (\u_j - u^\eps_j) \|_{L^2}^2 \leq C_\delta  \|\u\|_{\dot H^{\frac{3}{2}+\delta}} \left( \| \u-u^\eps \|_{\dot H^{\frac{1}{2}}} + \| \u-u^\eps \|_{\dot H^{\frac{3}{2}}}\right).
$$
Immediately, the remaining term estimates
\begin{equation}
 \sum_q 2^q \| \Delta_q \partial_i (\u_j - u^\eps_j) \, . \, S_q \u_i  \|_{L^2}^2 \leq C  \|\u-u^\eps\|_{\dot H^{\frac{3}{2}}}^2  \|\u\|_{L^\infty}^2.
\end{equation}
Finally, we obtain
\begin{eqnarray*}
 \eps \| A \|_{L^2_T \dot H^{\frac{1}{2}}}^2 & \! \leq  &\! C_{\delta , \eta} \eps \int_0^T \!\!\! \left( \|\u\|_{\dot H^{\frac{3}{2}+\delta}}^2 + \|u^\eps\|_{\dot H^{\frac{3}{2}+\delta}}^2 \right) \, \mathcal{E}_{\eps , \alpha ,u^\eps}(t) \, \d t + \\
 & \! + & \! C \eps \left( \| \u \|_{L^\infty}^2 + \|u^\eps\|_{L^\infty}^2 \right) \|\u-u^\eps \|_{L^2_T \dot H^{\frac{3}{2}}}^2 + \\
 & \! + & \! C_\delta  \eta \eps \left( \| \u \|_{L^\infty_T \dot H^{\frac{3}{2}+\delta}}^2 +  \|u^\eps\|_{L^\infty_T \dot H^{\frac{3}{2}+\delta}}^2 \right)  \| \u-u^\eps \|_{L^2_T \dot H^{\frac{3}{2}}}^2 ,
\end{eqnarray*}
where we use that $\eps \| \u \|_{L^\infty_T \dot H^{\frac{3}{2}+\delta}}^2 = \mathcal{O}(1)$ (see subsection \ref{rappels3d}).

\subsubsection{Estimate on $\int_0^T \int_{\R^3} \Lambda^{\frac{1}{2}} \nabla p^\eps . \Lambda^{\frac{1}{2}} \u \, \d x \, \d t $ }
Applying the $\div$ operator to $(HNS^\eps)$, we obtain the identity
$$ p^\eps = - \frac{1}{\Delta} \div (u^\eps. \nabla)u^\eps .  $$
Then, recall that $\div \u \in L^2_T \dot H^{\frac{1}{2}}$ and $\| \div \u \|_{L^2_T \dot H^{\frac{1}{2}}} = \mathcal{O}(\sqrt{\alpha})$. Now, we easily estimate the integral as follows:
\begin{eqnarray*}
 \left| \int_0^T \!\!\!\int_{\R^3} \Lambda^{\frac{1}{2}} \nabla p^\eps . \Lambda^{\frac{1}{2}} \u \, \d x \, \d t \right|  & = & \left| \int_0^T \!\!\!\int_{\R^3} \Lambda p^\eps . \div \u \, \d x \, \d t \right| \\
 & \leq & \| \Lambda p^\eps \|_{L^2_T L^{\frac{3}{2}}} \, \| \div \u \|_{L^2_T L^3} \\
 & \leq & C \sum_{i,j,k=1}^3  \| \big( \partial_k u^\eps_i \big) u^\eps_j \|_{L^2_T L^{\frac{3}{2}}} \, \| \div \u \|_{L^2_T L^3}  \\
 & \leq & C \sum_{i,j,k=1}^3 \| \partial_k u^\eps_i\|_{L^2_T L^3} \, \|u^\eps_j\|_{L^\infty_T L^3} \, \| \div \u \|_{L^2_T L^3} \\
 & \leq & C  \| u^\eps \|_{L^2_T \dot H^{\frac{3}{2}}} \, \|u^\eps\|_{L^\infty_T \dot H^{\frac{1}{2}}} \, \| \div \u \|_{L^2_T \dot H^{\frac{1}{2}}} \\
% & \leq & \frac{C}{\sqrt{\alpha}} \| \div \u \|_{L^2_T \dot H^{\frac{1}{2}}}^2  \,+ \, \frac{\sqrt{\alpha}}{2}  \| u^\eps \|_{L^2_T \dot H^{\frac{3}{2}}}^2 \, \|u^\eps\|_{L^\infty_T \dot H^{\frac{1}{2}}}^2 
 & = & \mathcal{O}(\sqrt{\alpha}).
\end{eqnarray*}

\subsubsection{Estimate on $ 2 \eps \int_0^T \int_{\R^3} \Lambda^{\frac{1}{2}} \nabla p^\eps . \Lambda^{\frac{1}{2}} \partial_t \u \, \d x \, \d t $}
First, two integrations by parts (one in space and another in time) give
$$\MakeUppercase{\romannumeral 1}:=2 \eps \int_0^T \int_{\R^3} \Lambda^{\frac{1}{2}} \nabla p^\eps . \Lambda^{\frac{1}{2}} \partial_t \u \, \d x \, \d t =  2\eps \int_0^T \int_{\R^3} \Lambda^{\frac{1}{2}} \partial_t p^\eps . \Lambda^{\frac{1}{2}} \div \u  \, \d x \, \d t $$
which is easier to estimate. \\
Notice that $\Lambda^{\frac{1}{2}} \div \u \in L^2_T L^2(\R^3)$.
Now, recall that 
$$ \Lambda^{\frac{1}{2}} \partial_t p^\eps = \Lambda^{\frac{1}{2}} \partial_t  \frac{1}{\Delta} \div (u^\eps.\nabla)u^\eps.$$ 
So we have
\[ \| \Lambda^{\frac{1}{2}} \partial_t p^\eps \|_{L^2_T L^2} \leq C  \| \partial_t (u^\eps.\nabla)u^\eps \|_{L^2_T \dot H^{-\frac{1}{2}}} = C \| \partial_t \big( u^\eps \otimes u^\eps \big) \|_{L^2_T \dot H^{\frac{1}{2}}} \leq C\| \partial_t u^\eps \otimes u^\eps \|_{L^2_T \dot H^{\frac{1}{2}}}.
\]
 We shall estimate this term by a dyadic Littlewood-Paley decomposition. Using Lemma \ref{paraproduct} in the appendix, we have
$$ \| u^\eps \otimes \partial_t u^\eps \|_{\dot H^{\frac{1}{2}}}^2 \leq \sum_{i,j} \sum_{p  \in \mathbb{Z}} 2^p \| \Delta_p u^\eps_i \, S_{p+1} \partial_t u^\eps_j \|_{L^2}^2 + \sum_{i,j} \sum_{q \in \mathbb{Z}} 2^q \| \Delta_q \partial_t u^\eps_j \, S_{q} u^\eps_i \|_{L^2}^2 . $$
Then we estimate each term separately, using Bernstein inequalities. First, we have
\begin{eqnarray*}
 \sum_p 2^p \| \Delta_p u^\eps_i \, S_{p+1} \partial_t u^\eps_j \|_{L^2}^2 &\leq & \sum_p 2^p \| \Delta_p u^\eps_i \|_{L^2}^2 \| S_{p+1} \partial_t u^\eps_j \|_{L^\infty}^2 \\
 & \leq & C \|u^\eps\|_{\dot H^{\frac{3}{2}}}^2  \sum_p 2^{-2p} \left\| \sum_{-1 \leq k \leq p} \Delta_k \partial_t u^\eps_j \right\|_{L^\infty}^2 \\
 & \leq & C \|u^\eps\|_{\dot H^{\frac{3}{2}}}^2  \sum_p 2^{-2p}  \sum_{-1 \leq k \leq p} 2^{3k} \left\| \Delta_k \partial_t u^\eps_j \right\|_{L^2}^2 \\
 & \leq & C  \|u^\eps\|_{\dot H^{\frac{3}{2}}}^2  \left\| \partial_t u^\eps \right\|_{\dot H^{\frac{1}{2}}}^2.
\end{eqnarray*}
Similarly, we obtain that
\[
 \sum_q 2^q \| \Delta_q \partial_t u^\eps_j \, S_{q} u^\eps_i \|_{L^2}^2  \leq  \sum_q 2^q \| \Delta_q \partial_t u^\eps_j \|_{L^2}^2   \| S_{q} u^\eps_i \|_{L^\infty}^2 
  \leq  C  \| \partial_t u^\eps \|_{\dot H^{\frac{1}{2}+\delta}}^2  \|u^\eps\|_{\dot H^{\frac{3}{2}+\delta}}^2.
\]
So finally, we have
\begin{eqnarray*}
 \MakeUppercase{\romannumeral 1} & \leq & C \eps \left( \|u^\eps\|_{L^\infty_T \dot H^{\frac{3}{2}}} \left\| \partial_t u^\eps_j \right\|_{L^2_T \dot H^{\frac{1}{2}}} + \| \partial_t u^\eps \|_{L^\infty_T \dot H^{\frac{1}{2}+\delta}} \|u^\eps\|_{L^2_T \dot H^{\frac{3}{2}+\delta}} \right) \|\div \u \|_{L^2_T \dot H^{\frac{1}{2}}} \\
 & = & \mathcal{O}\left( \sqrt{\alpha} \right).
\end{eqnarray*}

\subsubsection{Estimate on $\int_0^T \int_{\R^3} \Lambda (\u-u^\eps) . \left( u^\eps.\nabla u^\eps - \u.\nabla \u \right) \, \d x \, \d t $}
Let us write $u^\eps.\nabla u^\eps - \u.\nabla \u = (u^\eps-\u).\nabla u^\eps + \u.\nabla (u^\eps-\u)$. First, using the Sobolev embeddings $\dot H^1 \subset L^6 (\R^3)$ and $\dot H^{\frac{1}{2}} \subset L^3 (\R^3)$, we estimate the integral
$$ \tilde A := \int_0^T \!\!\!\!\! \int_{\R^3} \!\!\! \Lambda (\u-u^\eps) . (u^\eps-\u).\nabla u^\eps \, \d x \, \d t 
$$ as follows:
\begin{eqnarray*}
 \tilde A & \leq & \int_0^T \|\Lambda (\u-u^\eps) \|_{L^2} \, \|\u-u^\eps\|_{L^6} \, \|\nabla u^\eps\|_{L^3} \, \d t\\
 & \leq & C \int_0^T \|\u-u^\eps\|_{\dot H^1}^2 \, \|u^\eps\|_{\dot H^{\frac{3}{2}}} \, \d t \\
 & \leq & C \int_0^T \|\u-u^\eps\|_{\dot H^{\frac{1}{2}}} \, \|\u-u^\eps\|_{\dot H^{\frac{3}{2}}} \, \|u^\eps\|_{\dot H^{\frac{3}{2}}} \, \d t \\
 & \leq & \frac{C}{\eta} \int_0^T \|u^\eps\|_{\dot H^{\frac{3}{2}}}^2 \, \mathcal{E}_{\eps , \alpha ,u^\eps}(t) \, \d t + \eta \|\u-u^\eps\|_{L^2_T \dot H^{\frac{3}{2}}}^2
\end{eqnarray*}
Now, we are left with the term:
\begin{eqnarray*}
\Lambda (\u-u^\eps) .\u.\nabla (u^\eps-\u) &=& \sum_{i,j=1}^3 \partial_j (\u_i -u^\eps_i) \, (\u_i - u^\eps_i) \, \partial_i (\u_j - u^\eps_j) + \\
& &  + \sum_{i,j=1}^3 \partial_j (\u_i -u^\eps_i) \, u^\eps_i \, \partial_i (\u_j - u^\eps_j) = I_1 + I_2.
\end{eqnarray*}
One can easily check that the first part estimates as follows:
\begin{eqnarray*} 
\int_{\R^3} I_1 \, \d x & =&  \frac{1}{2} \sum_{i=1}^3 \int_{\R^3} \!\!\!\partial_i (\u_i -u^\eps_i)^2  \div \u \, \d x \\
& \leq & C \sum_{i=1}^3 \| \u_i -u^\eps_i \|_{L^3} \| \partial_i ( \u_i -u^\eps_i ) \|_{L^3} \| \div \u \|_{L^3} \\
& \leq & C \| \u -u^\eps \|_{\dot H^{\frac{1}{2}}} \| \u -u^\eps \|_{\dot H^{\frac{3}{2}}} \| \div \u \|_{\dot H^{\frac{1}{2}}} \\
& \leq & \frac{C \alpha}{\eta}  \mathcal{E}_{\epsalpha , u^\eps} + \eta \| \u - u^\eps \|_{\dot H^{\frac{3}{2}}}^2 ,
\end{eqnarray*}
where $\eta$ is small. Besides, by integrations by parts, we obtain for the second term
\begin{eqnarray*}
 \int_{\R^3} I_2 \, \d x & = & - \sum_{i,j=1}^3  \int_{\R^3} (\u_i -u^\eps_i)  \partial_j u^\eps_i  \partial_i (\u_j - u^\eps_j)   \, \d x - \sum_{i=1}^3 \int_{\R^3}  (\u_i - u^\eps_i) u_i \partial_i \div \u  \, \d x \\
 & = & - \sum_{i,j=1}^3  \int_{\R^3} (\u_i -u^\eps_i)  \partial_j u^\eps_i  \partial_i (\u_j - u^\eps_j)   \, \d x \\
 & & + \sum_{i=1}^3 \int_{\R^3}  \partial_i (\u_i - u^\eps_i) \, u^\eps_i \, \div \u   \, \d x \\
 & & + \sum_{i=1}^3 \int_{\R^3} (\u_i -u^\eps_i) \, \partial_i u^\eps_i \, \div \u  \, \d x \\
 & = & \MakeUppercase{\romannumeral 1} + \MakeUppercase{\romannumeral 2} + \MakeUppercase{\romannumeral 3} 
\end{eqnarray*}
We estimate these three terms as follows:
\begin{eqnarray*}
 \MakeUppercase{\romannumeral 1} & \leq & \sum_{i,j=1}^3 \|\u_i -u^\eps_i \|_{L^3} \, \|\partial_j u^\eps_i \|_{L^3} \, \|\partial_i (\u_j -u^\eps_j) \|_{L^3} \\
 & \leq & C \|\u-u^\eps\|_{\dot H^{\frac{1}{2}}} \, \|u^\eps\|_{\dot H^{\frac{3}{2}}} \, \|\u-u^\eps\|_{\dot H^{\frac{3}{2}}} \\
 & \leq & \frac{C}{\eta} \|u^\eps\|_{\dot H^{\frac{3}{2}}}^2 \, \mathcal{E}_{\eps , \alpha , u^\eps} \,+\, \eta \|\u-u^\eps\|_{\dot H^{\frac{3}{2}}}^2 .
\end{eqnarray*}

\begin{eqnarray*}
 \MakeUppercase{\romannumeral 2} & \leq & \sum_{i=1}^3 \|\partial_i (\u_i -u^\eps_i) \|_{L^3} \, \|u^\eps_i \|_{L^3} \, \|\div \u \|_{L^3} \\
 & \leq & C \|\u-u^\eps\|_{\dot H^{\frac{3}{2}}} \, \|u^\eps\|_{\dot H^{\frac{1}{2}}} \, \|\div \u \|_{\dot H^{\frac{1}{2}}} \\
 & \leq & C \|u^\eps_0\|_{\dot H^{\frac{1}{2}}}  \, \|\u-u^\eps\|_{\dot H^{\frac{3}{2}}} \, \|\div \u \|_{\dot H^{\frac{1}{2}}} \\
 & \leq & \eta \|\u-u^\eps\|_{\dot H^{\frac{3}{2}}}^2 + \frac{C}{\eta} \| \div \u\|_{\dot H^{\frac{1}{2}}}^2 .
\end{eqnarray*}

\begin{eqnarray*}
 \MakeUppercase{\romannumeral 3} &\leq & \sum_{i=1}^3 \|\u_i -u^\eps_i \|_{L^3} \, \| \partial_i u^\eps_i \|_{L^3} \, \|\div \u\|_{L^3} \\
 & \leq & C \|\u-u^\eps\|_{\dot H^{\frac{1}{2}}}^2 \, \|u^\eps\|_{\dot H^{\frac{3}{2}}}^2 \, + \, \frac{1}{2} \|\div \u \|_{\dot H^{\frac{1}{2}}}^2 \\
 & \leq & C \|u^\eps\|_{\dot H^{\frac{3}{2}}}^2 \, \mathcal{E}_{\eps , \alpha , u^\eps} \,+\,  \frac{1}{2} \|\div \u \|_{\dot H^{\frac{1}{2}}}^2 .
\end{eqnarray*}
Summarizing, we have obtained that
\begin{eqnarray*}
\int_0^T \!\!\!\!\! \int_{\R^3}\!\!\!\! \Lambda (\u-u^\eps) . \left( u^\eps.\nabla u^\eps - \u.\nabla \u \right) \, \d x \, \d t & \!\!\leq & \!\!\!\! \int_0^T \!\!\! \left(\! C_\eta  \|u^\eps\|_{\dot H^{\frac{3}{2}}}^2 + C_\eta \alpha \! \right)  \mathcal{E}_{\eps , \alpha , u^\eps}(t)  \d t +  \\
& & + C\eta \|\u-u^\eps\|_{L^2_T \dot H^{\frac{3}{2}} }^2  +  \mathcal{O}(\alpha) .
\end{eqnarray*}

\subsubsection{Conclusion}
Since $\mathcal{E}_{\epsalpha , u^\eps}(0)=0$, if we choose $\eta$ and $\eps$ small enough, we obtain the estimate
$$
\mathcal{E}_{\epsalpha , u^\eps}(T) \leq \mathcal{O}(\sqrt{\alpha})  + \! \int_0^T \!\!\! \left[ C_{\delta , \eta} \eps \left( \| \u\|_{\dot H^{\frac{3}{2}+\delta}}^2 + \| u^\eps \|_{\dot H^{\frac{3}{2}+\delta}}^2 \right) + C_\eta \|u^\eps \|_{\dot H^{\frac{3}{2}}}^2 + C_\eta \alpha \right] \! \mathcal{E}_{\epsalpha , u^\eps}(t)  \d t .
$$
Now, notice that $u^\eps \in L^2_T \dot H^{\frac{3}{2}} \cap L^2_T \dot H^{\frac{3}{2}+\delta} $ and that $\eps \| \u \|^2_{L^2_T \dot H^{\frac{3}{2}+\delta}} = \mathcal{O}(1)$ then apply the Gronwall's lemma and obtain that, for all positive $T$,
$$\mathcal{E}_{\epsalpha , u^\eps}(T) = \mathcal{O}(\sqrt{\alpha}) .$$ 
As in the 2D case, Theorem \ref{th} concludes the proof of Theorem \ref{th_alpha}.

\section{Aknowledgements}
I would like to thank my supervisors Valeria Banica and Pierre-Gilles Lemarié-Rieusset for their help and advice.

\begin{appendices}
%\makeatletter
%\def\@seccntformat#1{Appendix~\csname the#1\endcsname:\quad}
%\makeatother
\section{Littlewood-Paley theory}  \label{lptheory}
One of the main tools we use in this paper is the dyadic Littlewood-Paley decomposition. In this subsection, we briefly recall some important results. Our main references for the subject are the books by Alinhac and Gérard \cite{alinhac} and by Lemarié-Rieusset \cite{pglr}. \\
\\
First, recall that the homogeneous Sobolev norm $\dot H^s$ writes
\[
 \|u\|_{\dot H^s}^2 = \int_{\R^d} |\xi|^{2s} |\hat u (\xi)|^2 \, \d \xi = \sum_{p \in \mathbb{Z}} \int_{2^p \leq |\xi | \leq 2^{p+1}} |\xi|^{2s} |\hat u (\xi)|^2 \, \d \xi.
\]
In particular, we can estimate this norm in terms of the $L^2$ norm of $u|_{\{ 2^p \leq |\xi | \leq 2^{p+1} \}}$ as follows:
\[
 \sum_{p \in \mathbb{Z}} 2^{2ps} \int_{2^p \leq |\xi | \leq 2^{p+1}}  |\hat u (\xi)|^2 \, \d \xi \leq \|u\|_{\dot H^s}^2 \leq 2^{2s} \sum_{p \in \mathbb{Z}} 2^{2ps} \int_{2^p \leq |\xi | \leq 2^{p+1}} |\hat u (\xi)|^2 \, \d \xi .
\]
Now, we shall approximate the functions $\mathbf{1}_{\{ 2^p \leq |\xi | \leq 2^{p+1} \}}$ by smooth ones. So let us consider a nonnegative function $\varphi \in C_0^\infty (\R^d)$ such that $\varphi (\xi) = 1$ if $|\xi|\leq 1$ and $\varphi(\xi)=0$ if $|\xi| \geq 1+\eps$. Then define the function
\[ \psi : \xi \mapsto \varphi \left( \frac{\xi}{2} \right) - \varphi (\xi) \]
which is supported in the annulus $\left\lbrace 1 \leq |\xi| \leq 2(1+\eps) \right\rbrace$.\\ 
%For all $u \in \mathcal{S}(\R^d)$, that is $u$ in the Schwartz class, 
For all tempered distribution $u$, that is $u \in \mathcal{S}'(\R^d)$, 
we define $\Delta_p u$ the $p^{\textrm{th}}$ dyadic block of $u$  by
$$  \widehat{\Delta_p u}(\xi) = \psi (2^{-p}\xi) \hat u (\xi) $$
and we denote the partial sums by
$$ S_p u = \sum_{q < p} \Delta_q u $$
for all $p \in \mathbb{Z}$.\\ %Every dyadic block is in $H^{+\infty}$ once $u \in H^s$.\\
%\begin{rk} The definition of $\Delta_p u$ holds also for $u\in L^2$ and for $\hat u \in L^1_{\textrm{loc}}$. \end{rk}
Now, notice that $ \forall \xi \neq 0 ~,~~  \sum_{p\in \mathbb{Z}} \psi (2^{-p}\xi) = 1 $. So, in particular for all the functions $u$ considered in this work, we have
$$  \forall \xi \neq 0 ~,~~ \hat u (\xi) =  \sum_{p\in \mathbb{Z}} \widehat{\Delta_p u} (\xi) .  $$
\begin{rk} Every dyadic block is in $H^{+\infty}$ once $u \in H^s$. \end{rk}
We recall here some important properties of this dyadic decomposition. The results are stated for homogeneous Sobolev spaces but hold also for classical Sobolev spaces.
\begin{itemize}
 \item (Almost-orthogonality) For all $u \in L^2$, 
 $$ \sum_{p\in \mathbb{Z}} \| \Delta_p u \|_{L^2}^2 \leq \|u\|_{L^2}^2 \leq 2 \sum_{p\in \mathbb{Z}} \| \Delta_p u \|_{L^2}^2 .$$
 \item There exists a constant $C$ such that for all $s\in \R$ and $p \in \mathbb{Z}$,
 \begin{equation} \label{equivalence_des_normes} \frac{1}{C} \sum_{p \in \mathbb{Z}} 2^{2ps} \|\Delta_p u\|_{L^2}^2 \leq \|u \|_{\dot H^s}^2 \leq C \sum_{p \in \mathbb{Z}} 2^{2ps} \|\Delta_p u\|_{L^2}^2. \end{equation}
 We will use this property several times in this work. 
 \item There exists a constant $C$ such that, for all $\alpha \in \mathbb{N}^d$ and $p \in \mathbb{Z}$,
 $$\begin{array}{lcl}
     \| \partial^\alpha \Delta_p u \|_{L^2} ~ \leq  C 2^{p|\alpha|} \|u \|_{L^2} &,&  \| \partial^\alpha S_p u \|_{L^2} ~ \leq  C 2^{p|\alpha|} \|u \|_{L^2} ,\\
  \| \partial^\alpha \Delta_p u \|_{L^\infty}  \leq  C 2^{p|\alpha|} \|u \|_{L^\infty} &,&  \| \partial^\alpha S_p u \|_{L^\infty}  \leq  C 2^{p|\alpha|} \|u \|_{L^\infty} .
 \end{array}
 $$
\end{itemize}
Using this theory, we can prove some important product estimates. The following lemma tells how to write a product of tempered distributions:
\begin{lemme}[Paraproduct] \label{paraproduct}
 Given two tempered distributions $u,v \in \mathcal{S}'(\R^d)$, we can write their product, if it exists, as follows:
 \[
  uv = \sum_{p \in \mathbb{Z}} \Delta_p u \, S_{p+1} v + \sum_{q \in \mathbb{Z}} \Delta_q v \, S_q u .
 \]
\end{lemme}
Then, the following proposition is very useful to estimate some products.
\begin{prop} \label{tame}
 Let $s>0$ and $u,v \in L^\infty \cap \dot H^s$. Then the product $uv$ is also in $L^\infty \cap \dot H^s$ and
 $$ \|uv\|_{L^\infty \cap \dot H^s} \leq C \left( \|u\|_{L^\infty} \|v\|_{\dot H^s} + \|u\|_{\dot H^s} \|v\|_{L^\infty} \right).$$
\end{prop}
For a detailed proof of this proposition, see the book by Alinhac and Gérard.
%\begin{lemme}
%  Let $(a_q)_{q\geq -1}$ be a sequence of functions such that 
%  $$ \left\lbrace \begin{array}{l}
%  \hat a_q ~ \textrm{is supported in} ~ \left\lbrace |\xi| \leq C \, 2^q \right\rbrace \\
%  \exists s>0 , ~ \exists (c_q) \in \ell^2 ,~ \sum_{q\geq -1} c_q^2 \leq 1 ~~:~~ \|a_q\|_{L^2} \leq \tilde C \, c_q \, 2^{-qs}
%  \end{array}
%  \right.
%  $$
% Then $u= \sum_{q\geq -1} a_q \in H^s$ and $$\|u\|_{H^s} \leq C \tilde C .$$
% \end{lemme}
In this paper, we also use Bernstein inequalities:
\begin{prop}
 There exists a positive constant $C_0$ such that, for all $u \in \mathcal{S}'(\R^d)$, we have
 $$ \left\lbrace 
 \begin{array}{l}
 \forall k \geq 0 ,~ \forall p \geq 1 , ~~2^{qk} C_0^{-k} \|\Delta_q u \|_{L^p} \leq \sup_{|\alpha|=k} \| \partial^k \Delta_q u \|_{L^p} \leq 2^{qk} C_0^k \|\Delta_q u\|_{L^p} \\

 \forall p' \geq p \geq 1 , ~~ \|\Delta_q u\|_{L^{p'}} \leq C_0 2^{dq \left( \frac{1}{p} - \frac{1}{p'} \right)} \|\Delta_q u \|_{L^p}
 \end{array}
\right.
 $$
 and, for the partial sums:
 $$ \left\lbrace 
 \begin{array}{l}
 \forall k \geq 0 ,~ \forall p \geq 1 , ~~\sup_{|\alpha|=k} \| \partial^k S_q u \|_{L^p} \leq 2^{qk} C_0^k \|S_q u\|_{L^p} \\
 \forall p' \geq p \geq 1 , ~~ \|S_q u\|_{L^{p'}} \leq C_0 2^{dq \left( \frac{1}{p} - \frac{1}{p'} \right)} \|S_q u \|_{L^p}
 \end{array}
\right.
 $$
\end{prop}
Using the dyadic blocks, we can easily define the homogeneous Besov spaces $\dot{B}^s_{p,r}(\R^d)$, where $s\in \R$, $p,r\geq 1$ and $d\geq 1$,
$$ \dot{B}^s_{p,r}(\R^d) = \left\lbrace u \in \mathcal{S}'(\R^d) ~:~ \|u\|_{\dot B^s_{p,r}} := \left( \sum_{j\geq -1} 2^{jsr} \|\Delta_j u \|_{L^p}^r \right)^{\frac{1}{r}} \right\rbrace .$$
\end{appendices}

\bibliographystyle{alpha}
%\bibliography{article18052012}

\end{document}